\documentclass[11pt]{amsart}

\usepackage{amsthm, amsmath, amssymb, url, verbatim, xcolor, amscd}

\usepackage[normalem]{ulem}

\usepackage{hyperref}

\usepackage{latexsym}
\usepackage{amssymb}
\usepackage{amsmath}
\usepackage{amsfonts}
\usepackage[mathscr]{eucal}
\usepackage{color}  
\usepackage{tabulary}
\usepackage{stmaryrd}

\usepackage[all,color]{xy}
\usepackage{tikz}
\usetikzlibrary{matrix}

\usepackage{longtable}
\usepackage{rotating}
\usepackage{array}
\usepackage{multicol}
\usepackage{multirow}
\usepackage{makecell}
\usepackage{mathtools}
  \usepackage{cellspace}
\usepackage{booktabs}

\usepackage{enumitem}

\theoremstyle{plain}
\newtheorem{theorem}{Theorem}[section]
\newtheorem{corollary}[theorem]{Corollary}

\newtheorem{proposition}[theorem]{Proposition}

\newtheorem{main}{Theorem}
\newtheorem{Corollary}[main]{Corollary}

\newtheorem{thm}[theorem]{Theorem}
\newtheorem{cor}[theorem]{Corollary}
\newtheorem{lem}[theorem]{Lemma}
\newtheorem{prop}[theorem]{Proposition}
\theoremstyle{definition}
\newtheorem{rem}[theorem]{Remark}
\theoremstyle{definition}

\theoremstyle{remark}

\numberwithin{equation}{section}

\hypersetup{
colorlinks,
allcolors=blue
}

\newcommand{\N}{\mathbb{N}}
\newcommand{\Z}{\mathbb{Z}}
\newcommand{\Q}{\mathbb{Q}}
\newcommand{\R}{\mathbb{R}}

\DeclareMathOperator{\SO}{SO}
\DeclareMathOperator{\Or}{O}
\DeclareMathOperator{\SU}{SU}

\DeclareMathOperator{\Spin}{Spin}
\DeclareMathOperator{\Sp}{Sp}
\newcommand{\Gtwo}{\mathrm{G}_2}
\newcommand{\Ffour}{\mathrm{F}_4}

\newcommand{\ox}{\otimes}

\DeclareMathOperator{\codim}{codim}

\newcommand{\id}{\mathrm{id}}

\newcommand{\pd}{\partial}

\DeclareMathOperator{\im}{Im}

\DeclareMathOperator{\Hom}{Hom}
\DeclareMathOperator{\Diff}{Diff}
\DeclareMathOperator{\Ext}{Ext}
\DeclareMathOperator{\Homeo}{Homeo}
\DeclareMathOperator{\Sq}{Sq}

\newcommand{\sph}{\mathbf{S}}

\renewcommand{\leq}{\leqslant}
\renewcommand{\geq}{\geqslant}

\def\wt{\widetilde}

\def\x{\times}
\def\ox{\otimes}

\def\In{\subseteq}

\def\<{\langle}
\def\>{\rangle}

\def\bpm{\begin{pmatrix}}
\def\epm{\end{pmatrix}}
\def\bvm{\begin{vmatrix}}
\def\evm{\end{vmatrix}}
\def\bsm{\left(\begin{smallmatrix}}
\def\esm{\end{smallmatrix}\right)}

\def\beq{\begin{equation}}
\def\eeq{\end{equation}}

\begin{document}

 \title{Rational spheres and double disk bundles}

\author[J.\ DeVito]{Jason DeVito}
\address[J.\ DeVito]{Department of Mathematics, University of Tennessee at Martin, USA.}
\email{jdevito1@ut.utm.edu}

\author[M.\ Kerin]{Martin Kerin}
\address[M.Kerin]{Department of Mathematical Sciences, Durham University, UK}
\email{martin.p.kerin@durham.ac.uk}

\begin{abstract} A manifold $M$ is said to be a double disk bundle if it can be decomposed as a union of two disk bundles glued together by a diffeomorphism of their boundaries.  We show that if $M^n$ is a closed simply connected $n$-manifold with $n$ even which is simultaneously a double disk bundle and a rational homology sphere, then $M$ must be homeomorphic to a sphere.  In addition, we show that in any dimension, a highly connected rational homology sphere which supports a double disk bundle structure, then its "middle" cohomlogy group must be cyclic.
\end{abstract}

\maketitle

\section{Introduction}
Suppose $M$ is a smooth, $n$-dimensional manifold which admits a metric of non-negative sectional curvature.  Then the cohomology ring is ``small" in some sense.  One way this is made precise is through Cheeger and Gromoll's Soul Theorem \cite{CG}, which asserts that if $M$ is non-compact, it must have the structure of a vector bundle over a closed non-negatively curved manifold.  Thus, as an immediate corollary, its cohomology ring is finitely generated.  Refining this further, Gromov \cite{Gromov} found a bound on the total Betti number $\sum_{i=0}^n b_i$ which depends only on the dimension of $M$.  If one focuses instead on homotopy groups, the Bott conjecture asserts that such a manifold should be rationally elliptic: $\dim (\oplus_{k=2}^\infty \pi_k(M)\ox \Q) < \infty$.

From the perspective of (co)homology, the simplest closed manifolds are spheres.  One could argue that the next simplest examples are the \emph{rational homology spheres}, simply connected closed manifolds with the rational homology of a sphere.  In this context, a stark difference emerges depending on parity: there are many known examples of odd-dimensional rational homology spheres admitting a metric of non-negative sectional curvature and comparatively few examples in even dimensions.  

More specifically, in odd dimensions, apart from the standard spheres, one has numerous examples of homogeneous spaces and biquotients \cite{KZ}, cohomogeneity-one manifolds \cite{De,GVZ,GZ}, and quotients of these examples \cite{GKS1}.  In some dimensions, one has infinitely many examples distinct up to homotopy \cite{GKS1, GZ}, while the family $\SO(2n+1)/\SO(2n-1)$, $n \geq 1$, shows there exist examples in arbitrarily large dimensions (see, for example, \cite[Table 1]{KZ}).  Particularly interesting examples include the Wu manifold $\SU(3)/\SO(3)$ \cite{Ba, Cr}, the positively curved Berger space $\Sp(2)/\SU(2) = \SO(5)/\SO(3)$ \cite{Be, GKiS}, a positively curved exotic $T^1 \sph^4$ \cite{De, GVZ}, and all exotic $7$-spheres \cite{GKS1, GZ}.  

In contrast, in even dimensions, the only known rational homology spheres admitting metrics of non-negative sectional curvature are diffeomorphic to the standard spheres.

Recall that a manifold $M$ of dimension $n=2m+1$ is called \emph{highly connected} if $\pi_k(M) = 0$ for all $k < m$.  It follows from the Hurewicz theorem, Poincaré duality and the universal coefficient theorem that an odd-dimensional rational homology sphere is highly connected if and only if $H^k(M; \Z) =0$ for all $k \not\in \{0, m+1,2m+1\}$.  The group $H^{m+1}(M; \Z)$ is the obstruction to $M$ being homemorphic to a sphere and could, in principle, be any finite abelian group.  However, it turns out that $H^{m+1}(M;\Z)$ is cyclic in all known examples of highly connected, odd-dimensional rational homology spheres which admit non-negative sectional curvature.  For completeness, observe that, among the known examples of odd-dimensional rational homology spheres which admit non-negative sectional curvature, there are precisely three which are not highly connected:  the homogeneous spaces $\Gtwo/H$, where $H \in \{\SU(2)_3, \SO(3)_4, \SO(3)_{28}\}$ (using the notation from \cite{KZ}).  Indeed, from \cite{Mim} and the long exact sequence of homotopy groups for the fibrations $H \to \Gtwo \to \Gtwo/H$ it can be deduced that these $11$-dimensional rational homology spheres satisfy $\pi_2(\Gtwo/H) = \Z_2$ for $H \in  \{\SO(3)_4, \SO(3)_{28}\}$ and $\pi_2 (\Gtwo/\SU(2)_3) = 0$, $\pi_3(\Gtwo/\SU(2)_3) = \Z_3$, respectively.

The main goal of this article is to provide an explanation both for the lack of non-trivial examples of non-negatively curved, even-dimensional rational homology spheres and for the fact that the group $H^{m+1}(M^{2m+1}; \Z)$ is cyclic in all known examples of highly connected, odd-dimensional rational homology spheres which admit non-negative sectional curvature.   In a nutshell, the reason is that every known example of a non-negatively curved rational homology sphere has the structure of a \emph{double disk bundle} \cite{DGGK}; that is, up to diffeomorphism, the manifold can be obtained by gluing (the total spaces of) two disk bundles $DB_\pm \to B_\pm$ (possibly of different ranks) over closed manifolds $B_\pm$ together by a diffeomorphism of their boundaries.  For example, for each $1 \leq k \leq n$, the standard sphere $\sph^n$ can be obtained by gluing trivial disk bundles over $\sph^k$ and $\sph^{n-k-1}$ via the identity map.  

In fact, the \emph{Double Soul Conjecture} \cite{Gr} asserts that every closed, simply connected manifold $M$ with non-negative sectional curvature is a double disk bundle.  Evidence for the conjecture includes the Soul Theorem \cite{CG}, which can be viewed as a non-compact analogue, and the fact that all cohomogeneity-one manifolds admit this structure \cite{GGZ, Mostert}.  However, some caution is necessary: in the non-simply connected case, there exist infinitely many counterexamples \cite{DeV, Gr}, and, moreover, the conjecture is not even known to hold for compact Lie groups \cite{DGGK}.

While the main motivation for this work is geometric, the arguments are entirely topological.  It will be shown that requiring a rational homology sphere to be a double disk bundle imposes strong restrictions on the cohomology ring.  The consequence, in even dimensions, is that the only possible non-trivial examples are exotic spheres.

\begin{main}
\label{thm:main}
Let $M^{2m}$ be a  simply connected, even-dimensional, closed, smooth manifold with the rational cohomology of a sphere.  Then $M^{2m}$ admits a double disk-bundle decomposition if and only if it is homeomorphic (diffeomorphic for $m \in \{1, 2, 3\}$) to a sphere.
\end{main}

The case $m=2$ was previously shown by Ge and Radeschi \cite{GeRa} as a consequence of their classification of singular Riemannian foliations in dimension $4$.  The case $m=3$ was already shown by the authors in joint work with Galaz-Garc\'ia \cite{DGGK}.  For $m \geq 3$, the backward implication is immediate from \cite[Theorem H]{Sm2}.  For $m \geq 4$, the existence of exotic spheres in dimensions $\geq 7$ shows that ``homeomorphic'' cannot be replaced with ``diffeomorphic'' in general.

The proof of Theorem \ref{thm:main} leads to full homeomorphism classifications of the possible bases $B_\pm$ of the disk bundles and the boundaries $\partial DB_\pm$ along which the gluing is performed.  The following corollary says that, up to homeomorphism and with the exception of dimension four, the previously mentioned decomposition of a sphere into a union of trivial disk bundles over spheres is essentially all that can happen in the even-dimensional case.

\begin{Corollary}
\label{cor:homeotype}
Suppose $M^{2m}$ is an even-dimensional homotopy sphere that admits a double disk bundle decomposition $DB_- \cup_\partial DB_+$, where, without loss of generality, both of $B_\pm$ are assumed to be connected.  Then one of the following occurs:
\begin{enumerate}
\item  $B_\pm$ are points and $\partial DB_\pm$ are diffeomorphic to $\sph^{2m-1}$;

\item  $B_\pm$ are homotopy spheres and $\partial DB_\pm$ are homeomorphic to $B_- \times B_+$;

\item  $m=2$, $B_\pm$ are diffeomorphic to $\R P^2$ and $\partial DB_\pm$ are diffeomorphic to $\sph^3/Q_8$, where $Q_8 = \{\pm 1,\pm i, \pm j, \pm k\}$  is the quaternion group of order eight.
\end{enumerate}
\end{Corollary}

Turning attention now to the odd-dimensional case, recall that in dimension seven there exist infinitely many highly connected rational homology spheres which have the structure of a double disk bundle and admit non-negative sectional curvature \cite{GKS1, GZ}.  Each of these spaces has finite cyclic fourth cohomology group.  This is not an accident.

\begin{main}
\label{thm:main2}  
Let $M^{2m+1}$ be a highly connected, odd-dimensional, closed, smooth manifold with the rational cohomology of a sphere and suppose that $M^{2m+1}$ admits a double disk-bundle decomposition.  Then $H^{m+1} (M^{2m+1};\Z)$ is a (possibly trivial) finite cyclic group.
\end{main}

It is well known (see \cite{Bol, Mol}) that a simply connected, closed, smooth manifold which admits a singular Riemannian foliation (or, more generally, a transnormal system) of codimension one must be a double disk bundle.  In fact, the same holds if one replaces ``singular Riemannian foliation of codimension one''  with ``isoparametric foliation'' or ``manifold submetry onto a closed interval'' since, for simply connected, closed, smooth manifolds, these notions are all equivalent to the existence of a decomposition into the union of two \emph{linear} disk bundles glued along their common boundary (see \cite{KrLyRa, LY}).  Thus, the conclusions of Theorem \ref{thm:main}, Theorem \ref{thm:main2} and Corollary \ref{cor:homeotype} hold under any of these alternative hypotheses.  Observe, however, that Theorem \ref{thm:main}, Theorem \ref{thm:main2} and Corollary \ref{cor:homeotype} hold in more generality, as it is not assumed that the disk bundles are linear.

\

To summarize the strategy behind the proofs of Theorem \ref{thm:main} and Theorem \ref{thm:main2}, consider a closed, smooth manifold $M$ which admits a double disk-bundle decomposition $DB_- \cup_\partial DB_+$, where the gluing is along the common boundary $L \cong \partial DB_\pm$.  Grove and Halperin \cite{GH} determined the rational homotopy type and integral (co)homology groups of the homotopy fiber $F$ of the inclusion map $L\rightarrow  M$.  Together with the simple rational topology of both $M$ and the sphere bundles $L\rightarrow B_\pm$, this imposes very strong constraints on the rational topology of $B_\pm$ and $L$.

When $M$ is even dimensional, the main step, contained in Proposition \ref{P:one_sphere}, is to use the Leray-Hirsch Theorem to gain control of the integral cohomology ring of one of $B_\pm$, say $B_+$, and show that it must be a homotopy sphere.  The integral cohomology of $L$ is then easily determined from the sphere bundle $L\rightarrow B_+$.  Next, calculations involving Stiefel-Whitney classes and Steenrod operations are used to show that $B_-$ is either a homotopy sphere (in which case Theorem \ref{thm:main} follows easily from the Mayer-Vietoris sequence) or diffeomorphic to the Wu manifold $\SU(3)/\SO(3)$.  In the latter case, this forces $L$ to simultaneously be an $\sph^2$-bundle over $\SU(3)/\SO(3)$ and an $\sph^5$-bundle over $\sph^2$.  Lemma \ref{L:noWu} provides an argument using higher homotopy groups to show that this is impossible.

When $M$ is $(2m+1)$-dimensional, the highly connected assumption and knowledge of the integral cohomology of $F$ together yield great control over the integral cohomology of $L$.  In particular, via the Serre spectral sequence for the homotopy fibration $F\rightarrow L\rightarrow M$, it now follows easily that $H^{m}(M; \Z)$ is generated by at most two elements.  At this point, a detailed analysis of the commutative braid \eqref{E:braid} of exact sequences associated to the decomposition $DB_- \cup_L DB_+$ allows one to extract enough information to show that $H^m(M;\Z)$ is actually generated by at most one element.

The outline of the paper is as follows.  Section \ref{sec:bas} covers the necessary background, including information on characteristic classes of possibly non-linear sphere bundles, as well as information on the rational homotopy of double disk bundles.  Theorem \ref{thm:main} and Corollary \ref{cor:homeotype} are proved in Section \ref{S:evendim}, while Section \ref{S:highly_connected} is devoted to the proof of Theorem \ref{thm:main2}.

\

\textbf{Acknowledgments}:  The first author would like to thank the University of Oklahoma, where he was visiting when some of these results were obtained.  In addition, we would like to thank Jeffrey Carlson for comments on an earlier version of this article.  The first author was also supported by the NSF via both DMS-2105556 and DMS-2405266.  He is grateful for the support.  The second author is grateful to Marco Radeschi for conversations regarding singular Riemannian foliations and isoparametric foliations.


\section{Basic structure}\label{sec:bas}

In this section, we set up the basic terminology and structure necessary for the rest of the paper.  Throughout, the symbol $\cong$ will be used to indicate either that two manifolds are diffeomorphic or that two groups are isomorphic, depending on the context.  Homology and cohomology will be taken with integral coefficients, unless explicitly indicated otherwise.

\subsection{Double disk bundles}

Let $B_\pm$ be smooth, closed manifolds and suppose that $D^{\ell_\pm + 1} \to DB_\pm \to B_\pm$ are smooth disk bundles of rank $\ell_\pm + 1$ over $B_\pm$, respectively.  If there is a diffeomorphism $f \colon \partial DB_- \to \partial DB_+$ of the boundaries, then identifying $\partial DB_\pm$ via $f$ yields a smooth, closed manifold $DB_- \cup_f DB_+$, called a \emph{double disk bundle}.  Recall that smooth disk bundles admit smooth sections \cite[Section 6.7 and Theorem 12.2]{Steenrod}, so $B_\pm$ can always be considered as smooth submanifolds of $DB_- \cup_f DB_+$.  Moreover, by Whitehead's Theorem, the bundle projections induce homotopy equivalences $DB_\pm \simeq B_\pm$.

If $L$ denotes the common image of $\partial DB_\pm$ in $DB_- \cup_f DB_+$, then it can be viewed as the total space of two smooth sphere bundles $\sph^{\ell_\pm} \to L \to B_\pm$.  Indeed, the homomorphisms $\Diff(D^{\ell_\pm + 1}) \to \Diff(\sph^{\ell_\pm})$ given by restricting diffeomorphisms of $D^{\ell_\pm + 1}$ to the boundaries, yield descriptions of $\partial DB_\pm$ as the total spaces of the smooth $\sph^{\ell_\pm}$-bundles associated to the principal $\Diff(D^{\ell_\pm + 1})$-bundles underlying the disk bundles $DB_\pm$, respectively.  In other words, the space $L$ is the total space of two different smooth sphere bundles, whose structure groups reduce to $\Diff(D^{\ell_\pm + 1})$, respectively.

An arbitrary smooth, closed, connected manifold $M$ is said to admit a \emph{double disk-bundle decomposition} if there exists a diffeomorphism $\Phi \colon M \to DB_- \cup_f DB_+$ from $M$ to a double disk bundle $DB_- \cup_f DB_+$.  By an abuse of notation, $B_\pm$ and $L$ will be used to denote the images of pulling back to $M$ via $\Phi$ the corresponding objects in the double disk bundle $DB_- \cup_f DB_+$, while the decomposition itself will typically be denoted by $DB_- \cup_L DB_+$, unless knowledge of the precise gluing map is needed.  As a consequence of \cite[Proposition 4.1]{DGGK}, it may always be assumed that $B_\pm$ are connected.

Recall that the diffeomorphism group $\Diff(\sph^k)$ deformation retracts onto $\Or(k+1)$ whenever $k \leq  3$ \cite{BK, Ha1, Sm1}.  Hence, it may be implicitly assumed that $\sph^{\ell_\pm} \to L \to B_\pm$ is a linear bundle if $\ell_\pm \leq 3$, respectively.  The inclusion $L \to M$ gives rise to an additional homotopy fibration $F \to L \to M$, where $F$ denotes the so-called \emph{homotopy fiber}.

We will borrow terminology from the related study of singular Riemannian foliations and often refer to $B_\pm$ and $L$ as the \emph{singular} and \emph{regular leaves}, respectively, of the double disk-bundle decomposition of $M$.


\subsection{The topology of double disk bundles} 

There is a commutative braid diagram \cite{Wall}
\begin{equation}
\label{E:braid}
\scalebox{0.78}{
\xymatrix@=0.4cm{
\phantom{X} \ar[dr]^(0.4){} \ar@(ur,ul)[rr] & 
& H^j(B_+) \ar[dr]^{} \ar@(ur,ul)[rr]^{} & 
& H^{j+1}(M,B_+) \ar[dr]^(0.55){} \ar@(ur,ul)[rr] & 
&  H^{j+1}(B_-) \ar[dr]^(0.4){} \ar@/^/@{-}[r]+/u 5mm/ & \\
& H^j(M)  \ar[ur]^{} \ar[dr]^(0.55){} &
& H^j(L) \ar[ur]^(0.45){} \ar[dr]^(0.5){} & 
& H^{j+1}(M) \ar[ur]^(0.55){} \ar[dr]^(0.55){} &  
& \phantom{X} \\
\phantom{X} \ar[ur]^(0.4){} \ar@(dr,dl)[rr] & 
& H^j(B_-) \ar[ur]^{} \ar@(dr,dl)[rr]_{}  & 
& H^{j+1}(M,B_-) \ar[ur]^(0.5){} \ar@(dr,dl)[rr] & 
& H^{j+1}(B_+) \ar[ur]^(0.5){} \ar@/_/@{-}[r]+/d 5mm/ & 
}}
\end{equation}
associated to each double disk-bundle decomposition $M^n \cong DB_- \cup_L DB_+$, where each strand is the long exact sequence of a pair and, in particular, the isomorphisms $H^j(M, B_\pm) \cong H^j(M,DB_\pm) \to H^j(DB_\mp,L)$ given by excision are used implicitly in the diagram.  Therefore, the homomorphisms $H^j(L) \to H^{j+1}(M, B_\mp)$ correspond to the boundary homomorphisms in the long exact sequences for the pairs $(DB_\pm, L)$.  Furthermore, if $M^n$ is orientable, then the disk bundles $DB_\pm \In M^n$ are $n$-dimensional, compact, orientable manifolds with boundary $\pd DB_\pm = L$.  Thus, Poincar\'e duality for manifolds with boundary (together with excision) yields 
\begin{equation}
\label{E:isoms}
H^{j}(M, B_\pm) \cong H^{j}(DB_\mp, L) \cong H_{n-j}(B_\mp) \,.
\end{equation}
Note that it is a simple exercise in diagram chasing to derive the Mayer-Vietoris sequence in cohomology for $M^n \cong DB_- \cup_L DB_+$ from the braid diagram \eqref{E:braid}.

It will be important to determine when a manifold admitting a double disk-bundle decomposition is simply connected.  The following proposition, a proof of which can be found in \cite[Proof of Proposition 3.5]{GH} or (using the language of  cohomogeneity-one manifolds) in \cite[Proposition 1.8]{Ho}, shows that this can be determined from the sphere bundles associated to the decomposition.  

\begin{proposition}
\label{prop:pi1}  
Suppose that a smooth, closed manifold $M$ admits a double disk-bundle decomposition $DB_- \cup_L DB_+$, where $B_\pm$ are both connected and the associated sphere bundles $\sph^{\ell_\pm}\rightarrow L \rightarrow B_\pm$ have $\ell_\pm \geq 1$.  Then $M$ is simply connected if and only if $\pi_1(L)$ is generated by the images of the induced homomorphisms $\pi_1(\sph^{\ell_\pm}) \to \pi_1(L)$ in the respective long exact homotopy sequences for the bundles.
\end{proposition}

\begin{rem}
\label{R:codim}
It follows from \cite[Proposition 4.2]{DGGK} that a double disk-bundle decomposition $DB_- \cup_L DB_+$ of a smooth, closed, simply connected manifold must have both $\ell_\pm \geq 1$ whenever  $B_\pm$ are connected.
\end{rem}

In \cite{GH}, Grove and Halperin studied spaces admitting double disk-bundle decompositions from the perspective of rational homotopy theory.  One of their main results can be summarized as follows:

\begin{theorem}[Grove--Halperin \protect{\cite{GH}}]
	\label{T:HOM_FIBER}
Suppose that a smooth, closed, simply connected manifold $M$ admits a double disk-bundle decomposition $DB_- \cup_L DB_+$, where $B_\pm$ are both connected and $L \cong \partial D B_\pm$.  If $F$ denotes the homotopy fiber of the inclusion $L \to M$, then $L$ and $F$ are nilpotent spaces and a finite cover of 
$F$ is rationally homotopy equivalent to one of the spaces listed in Table \ref{table:Qtype}, where $A_m(4)$ denotes a certain simply connected topological space whose non-trivial rational homotopy groups are in degrees $4$, $7$ and $4m - 1$.

Moreover, the possible fundamental groups and cohomology groups of $F$ are listed in Table \ref{table:Qtype} and Table \ref{table:cohomF}, respectively, arranged according to the codimensions of $B_\pm$ in $M$ and the orientability of the sphere bundles $\sph^{\ell_\pm} \to L \to B_\pm$ associated to the decomposition $DB_- \cup_L DB_+$.
\end{theorem}


\begin{table}
\centering
\resizebox{\columnwidth}{!}{%
\begin{tabular}{|Sc|Sc|Sc|Sc|}
\hline
\multirow{2}{*}{$\pi_1(F)$} &
\multirow{2}{*}{$F \simeq_\Q$} &
\multirow{2}{*}{$\{\alpha, \beta\} = \{\ell_\pm\}$} &
\multicolumn{1}{>{\centering\let\newline\\\arraybackslash\hspace{0pt}}p{35mm}|}{Orientability of \newline  $\sph^{\ell_\pm}$-bundles} \\
\hline \hline

$\Z^2$ & $\sph^1 \x \sph^1 \x \Omega \sph^3$ & \multirow{4}{*}{$1 = \alpha = \beta$} & Both \\ 
\cline{1-2} \cline{4-4}

$\Z \oplus \Z_2$ & $\sph^1 \x \sph^3 \x \Omega \sph^5$ &  & One \\
\cline{1-2} \cline{4-4}

$Q_8$ & $\sph^3 \x \sph^3 \x \Omega \sph^7$ &  & Neither \\
\hline

\multirow{2}{*}{$\Z$} & $\sph^1 \x \sph^\beta \x \Omega \sph^{\beta + 2}$ & $1 = \alpha < \beta$ & Both \\
\cline{2-4}
& $\sph^1 \x \sph^{2\beta + 1} \x \Omega \sph^{2\beta + 3}$ & $1 = \alpha < \beta$, $\beta$ odd & $\sph^1$-bundle \\
\hline

\multirow{11}{*}{$0$} & $\sph^\alpha \x \sph^\beta \x \Omega \sph^{\alpha + \beta + 1}$ & $1 < \alpha \leq \beta$ &  \multirow{11}{*}{Both} \\ \cline{2-3}
& $\sph^\alpha \x \Omega \sph^{\alpha + 1}$ & $1 < \alpha = \beta$  &  \\ \cline{2-3}
& $\SU(3)/T^2 \x \Omega \sph^7$ & \multirow{4}{*}{$2 = \alpha = \beta$} &  \\ \cline{2-2}
& $\Sp(2)/T^2 \x \Omega \sph^9$ & &  \\ \cline{2-2}
& $\Gtwo/T^2 \x \Omega \sph^{13}$ & &  \\ \cline{2-3}
& $\Sp(3)/\Sp(1)^3 \x \Omega \sph^{13}$ & \multirow{4}{*}{$4 = \alpha = \beta$} &  \\ \cline{2-2}
& $A_4(4) \x \Omega \sph^{17}$ & &  \\ \cline{2-2}
& $A_6(4) \x \Omega \sph^{25}$ & &  \\ \cline{2-3}
& $\Ffour/\Spin(8) \x \Omega \sph^{25}$ & $8 = \alpha = \beta$ &  \\ \hline
\end{tabular}
}
\vspace{10pt}
\caption{Properties of the homotopy fiber $F$ and the bundles $\sph^{\ell_\pm} \to L \to B_\pm$ associated to a double disk-bundle decomposition $DB_- \cup_L DB_+$}
\label{table:Qtype}

\end{table}


\begin{table}

\centering
\resizebox{\columnwidth}{!}{%
\begin{tabular}{|Sc|Sc|SlSl|}
\hline
\multirow{2}{*}{$\{\alpha, \beta\} = \{\ell_\pm\}$} &   
\multicolumn{1}{>{\centering\let\newline\\\arraybackslash\hspace{0pt}}p{35mm}|}{Orientability of \newline  $\sph^{\ell_\pm}$-bundles} &
\multicolumn{2}{c|}{\multirow{2}{*}{$H^j (F; \Z)$}}  
 \\
\hline \hline

\multirow{2}{*}{$\alpha \neq \beta$} & 
\multirow{2}{*}{Both} & 
$\Z,$ & 
$j = 0$, or $j \in \{\alpha, \beta\} \!\! \mod \alpha + \beta$ \\ 
& 
& 
$\Z^2,$ & 
$j > 0$ and $j \equiv 0 \!\! \mod \alpha + \beta$ \\ \hline

\multirow{2}{*}{$\alpha = \beta$} & 
\multirow{2}{*}{Both} & 
$\Z,$ & 
$j = 0$ \\
& 
& 
$\Z^2,$ & 
$j > 0$ and $j \equiv 0 \!\! \mod \alpha$ \\ \hline

\multirow{4}{*}{$1 = \alpha < \beta$} & 
\multirow{4}{*}{$\sph^1$-bundle} & 
$\Z,$ & 
$j = 0$, or $j \equiv \pm 1 \!\! \mod 2 \beta + 2 $ \\ 
& 
& 
$\Z^2,$ & 
$j > 0$ and $j \equiv 0 \!\! \mod 2 \beta + 2 $ \\ 
& 
& 
$\Z_2,$ & 
$j \in \{ \beta + 1, \beta + 2 \} \!\! \mod 2 \beta + 2 $ \\ 
\hline

\multirow{5}{*}{$1 = \alpha = \beta$} & 
\multirow{5}{*}{One} & 
$\Z,$ & 
$j = 0$, or $j \equiv 1 \!\! \mod 4 $ \\ 
& 
& 
$\Z_2,$ & 
$j \equiv 2 \!\! \mod 4 $ \\ 
& 
& 
$\Z \oplus \Z_2,$ & 
$j \equiv 3 \!\! \mod 4$ \\
& 
& 
$\Z^2,$ & 
$j > 0$ and $j \equiv 0 \!\! \mod 4 $ \\ \hline

\multirow{4}{*}{$1 = \alpha = \beta$} & 
\multirow{4}{*}{Neither} & 
$\Z,$ & 
$j = 0$ \\ 
& 
& 
$\Z^2,$ & 
$j > 0$ and $j \equiv 0 \!\! \mod 3 $ \\ 
& 
& 
$\Z_2^2,$ & 
$j \equiv 2 \!\! \mod 3 $ \\ 
\hline

\end{tabular}
}
\vspace{10pt}
\caption{Non-trivial cohomology groups of the homotopy fiber $F$ of the inclusion $L = \partial DB_\pm \to M$ associated to a double disk-bundle decomposition $DB_- \cup_L DB_+$ of $M$}
\label{table:cohomF}

\end{table}

Note that, in \cite{GH}, Grove and Halperin listed the homology groups of the homotopy fiber $F$, from which one obtains Table \ref{table:cohomF} by an application of the universal coefficient theorem.  Furthermore, as described in Table \ref{table:Qtype}, Grove and Halperin proved that the homotopy fiber $F$ is always rationally homotopy equivalent to a nilpotent space having exactly one factor given by the loop space of a sphere.  If $\Omega \sph^k$ is such a loop-space factor, then it has a unique non-trivial rational homotopy group of even degree associated with it: namely, in degree $2s = k-1$, if $k$ is odd, and in degree $2s = 2(k-1)$ if $k$ is even.   The following theorem was proven in \cite[Theorem 4.11]{DGGK} by the authors together with F.\ Galaz-García.

\begin{theorem} 
\label{T:nontriv}
Suppose that $M$ is a smooth, closed, simply connected manifold which admits a double disk-bundle decomposition $DB_- \cup_L DB_+$ with $B_\pm$ connected. Then, with $2s$ as above and with reference to the long exact sequence of rational homotopy groups associated to the homotopy fibration $F\to L \to M$, the connecting homomorphism $\partial\colon\pi_{2s+1}^\Q(M)\to \pi_{2s}^\Q(F)$ is non-trivial.
\end{theorem}

Theorem \ref{T:nontriv} is particularly useful when $M$ has a unique non-trivial rational homotopy group of odd degree.  To that end, recall that spheres are intrinsically formal in the sense of rational homotopy theory \cite{FH1}.  In particular, every simply connected manifold $M$ with the rational homology of a sphere is rationally homotopy equivalent to a sphere, so that $\pi_k(M)\otimes \mathbb{Q}$ is non-trivial in precisely one odd degree.

\begin{cor}
\label{C:dims}
Suppose that $M$ is a simply connected rational homology sphere of dimension $n$ which admits a double disk-bundle decomposition $DB_- \cup_L DB_+$ with $B_\pm$ connected.  Then the homotopy fiber $F$ of the inclusion $L \to M$ is rationally homotopy equivalent to a nilpotent space having exactly one loop-space factor $\Omega \sph^k$, where 
\begin{enumerate}
\item $k = n$,  or  
\item $k = 2n-1$, with $n$ even, or
\item $k = \frac{n+1}{2}$, with $n \equiv 3 \!\! \mod 4$.
\end{enumerate} 
\end{cor} 

\begin{proof}
Since $M$ is a rational homology sphere, the unique rational homotopy group $\pi_{2s+1}^\Q(M)$ for $M$ of odd degree occurs when $2s = n-1$, for $n$ odd, or when $2s = 2n-2$, for $n$ even.  By Theorem \ref{T:nontriv}, this implies that the loop-space factor $\Omega \sph^k$ of $F$ must satisfy one of the following conditions:
\begin{enumerate}[label=(\roman*)]
\item $n$, $k$ both odd and $k = 2s+1 = n$, or
\item $n$ odd, $k$ even and $k = s+1 = (n+1)/2$, or
\item $n$ even, $k$ odd and $k = 2s+1 = 2n - 1$, or
\item $n$, $k$ both even and $k = s + 1 = n$.
\end{enumerate}
The result now follows immediately.
\end{proof}

By combining Corollary \ref{C:dims} with Tables \ref{table:Qtype} and \ref{table:cohomF}, there are additional restrictions on the cohomology groups of the homotopy fiber $F$ whenever $M$ is a rational homology sphere.  Of particular interest later will be the situations where $F$ can have torsion in its cohomology.

\begin{cor}
\label{C:torF}
Suppose that $M$ is a simply connected rational homology sphere of dimension $n$ which admits a double disk-bundle decomposition $DB_- \cup_L DB_+$ with $B_\pm$ connected.  Then all of the cohomology groups of the homotopy fiber $F$ of the inclusion $L \to M$ are free abelian of rank $\leq 2$, unless one of the following holds:
\begin{enumerate}
\item 
\label{i:dim4}
$n = 4$ and $\codim(B_\pm) = 2$, where neither of the circle bundles $\sph^1 \to L \to B_\pm$ is orientable, in which case 
$$
H^j(F) \cong \begin{cases}
\Z, & j = 0\,, \\
\Z^2, & j > 0 \ \text{ and } \ j \equiv 0 \!\! \mod 3\,, \\
0, & j \equiv 1 \!\! \mod 3\,, \\
\Z_2^2, & j \equiv 2 \!\! \mod 3\,;
\end{cases}
$$

\item 
\label{i:dim5}
$n = 5$ and $\codim(B_\pm) = 2$, where exactly one of the circle bundles $\sph^1 \to L \to B_\pm$ is orientable, in which case 
$$
H^j(F) \cong \begin{cases}
\Z, & j = 0 \ \text{ or } \ j \equiv 1 \!\! \mod 4 \,, \\
\Z^2, & j > 0 \ \text{ and } \ j \equiv 0 \!\! \mod 4\,, \\
\Z_2, & j \equiv 2 \!\! \mod 4 \,, \\
\Z \oplus \Z_2, & j \equiv 3 \!\! \mod 4\,;
\end{cases}
$$

\item  
\label{i:dim7}
$n = 7$, and $\codim(B_\pm) = 2$, where neither of the circle bundles $\sph^1 \to L \to B_\pm$ is orientable, in which case
$$
H^j(F) \cong \begin{cases}
\Z, & j = 0\,, \\
\Z^2, & j > 0 \ \text{ and } \ j \equiv 0 \!\! \mod 3\,, \\
0, & j \equiv 1 \!\! \mod 3\,, \\
\Z_2^2, & j \equiv 2 \!\! \mod 3\,;
\end{cases}
$$

\item 
\label{i:dim2n+3}
$n = 2 \beta + 3 $, with $\beta > 1$ odd, and $\{\codim(B_\pm)\} = \{2, \beta+1\}$, where only the corresponding circle bundle is orientable, in which case
$$
H^j(F) \cong  \begin{cases}
\Z, & j = 0 \ \text{ or } \ j \equiv \pm 1 \!\! \mod n-1 \,, \\
\Z^2, & j > 0 \ \text{ and } \ j \equiv 0 \!\! \mod n-1 \,, \\
\Z_2, & j \in \{ \frac{n-1}{2}, \frac{n+1}{2} \} \!\! \mod n-1\,, \\
0, & \text{otherwise.}
\end{cases}
$$
\end{enumerate} 
\end{cor} 

\begin{proof}
From Tables \ref{table:Qtype} and \ref{table:cohomF}, it is clear that all cohomology groups of the homotopy fiber $F$ are free abelian groups of rank $\leq 2$ unless $F$ and the sphere bundles $\sph^{\ell_\pm} \to L \to B_\pm$ satisfy one of the following scenarios:
\begin{enumerate}[label=(\roman*)]
\item 
\label{i:k5}
$\ell_\pm = 1$, with only one of the bundles $\sph^1 \to L \to B_\pm$ being orientable, and $F \simeq_\Q \sph^1 \x \sph^3 \x \Omega \sph^5$; or 

\vspace*{1mm}
\item 
\label{i:k7}
$\ell_\pm = 1$, with neither of the bundles $\sph^1 \to L \to B_\pm$ being orientable, and $F \simeq_\Q \sph^3 \x \sph^3 \x \Omega \sph^7$; or 

\vspace*{1mm}
\item 
\label{i:k2b+3}
$\{\ell_\pm\} = \{1, \beta\}$, where $\beta > 1$ is odd and only the circle bundle among $\sph^{\ell_\pm} \to L \to B_\pm$ is orientable, and $F \simeq_\Q \sph^1 \x \sph^{2\beta + 1} \x \Omega \sph^{2\beta + 3}$. 
\end{enumerate}

Note that the loop-space factor $\Omega \sph^k$ of $F$ has $k$ odd in all three scenarios; namely, $k \in \{5, 7, 2 \beta + 3\}$ for some odd $\beta > 1$.  Thus, it follows from Corollary \ref{C:dims} that either $k = 2n-1$, with $n$ even, or $k = n = \dim (M)$.

If $2n-1 = k \in \{5, 7, 2 \beta + 3\}$, where $n$ is even and $\beta$ odd, then the only possibility is $n=4$ and $k = 7$, corresponding to case \ref{i:k7} above.  Otherwise, $n = k \in \{5, 7, 2 \beta + 3\}$, with all three cases \ref{i:k5}--\ref{i:k2b+3} possible.  The cohomology groups of the homotopy fiber $F$ can now be read off Table \ref{table:cohomF} easily in each of these four cases. 
\end{proof}

In the case of Corollary \ref{C:torF}\eqref{i:dim4}, Poincar\'e duality and the universal coefficient theorem force $M$ to be an integral homology $4$-sphere and, hence, homeomorphic to $\sph^4$, by Freedman's classification \cite{Fr}.  Now, since $M$ admits a double disk-bundle decomposition, it must, in fact, be diffeomorphic to $\sph^4$ \cite{GeRa}.  In particular, the scenario described in Corollary \ref{C:torF}\eqref{i:dim4} arises for a well-known cohomogeneity-one action on $\sph^4$ with codimension-two singular orbits; see \cite{Par} for details.  

With respect to Corollary \ref{C:torF}\eqref{i:dim5}, recall from \cite[Theorem B]{DGGK} that, if $M$ is a simply connected, $5$-dimensional rational homology sphere which admits a double disk-bundle decomposition, it must be diffeomorphic to either $\sph^5$ or the Wu manifold $\SU(3)/\SO(3)$.  In each case, there exists a cohomogeneity-one action with codimension-two singular orbits such that Corollary \ref{C:torF}\eqref{i:dim5} applies; see \cite{Ho} for details.

The scenario in Corollary \ref{C:torF}\eqref{i:dim7} arises for each member of the infinite family of rational $7$-spheres constructed in \cite{GKS1}.  This includes all exotic spheres in dimension seven.  Finally, the case of Corollary \ref{C:torF}\eqref{i:dim2n+3} arises via the tensor product action of $\Or(\beta+2)\times \Or(2)$ on $\sph^{2\beta+3}\subseteq \mathbb{R}^{2\beta+4}\cong \mathbb{R}^{\beta+2}\otimes \mathbb{R}^2$ \cite[(12.8)]{Wang}.


\subsection{Characteristic classes of smooth bundles}
\label{S:nonlin}

This section gives a brief overview of the extension of the usual theory of characteristic classes (typically defined for vector bundles, linear sphere bundles, or linear disk bundles) to the more general class of smooth bundles.  Although the discussion focuses on smooth bundles, it is perhaps worth noting that the definitions and properties of the Euler class and the Stiefel-Whitney classes carry over verbatim to topological disk and sphere bundles.

Let $D^{k+1} \to E \xrightarrow{\pi} B$ be a smooth disk bundle over a closed, connected manifold $B$.  Let $\partial E$ be the boundary of the total space $E$ and let $\sph^k \to \partial E \xrightarrow{\pi} B$ be the associated smooth sphere bundle.  Finally, if the bundle $D^{k+1} \to E \xrightarrow{\pi} B$ is orientable, let $R$ be any commutative ring with identity; otherwise, let $R = \Z_2$.  By the Thom Isomorphism Theorem \cite{Thom} (see also \cite{Co} and \cite[Section 4.D]{Hatcher}), there exists a Thom class $u_R \in H^{k+1}(E, \partial E; R)$ such that the homomorphism
\begin{align*}
\Phi : H^j(B; R) &\to H^{j + k + 1}(E, \partial E; R) \\ 
x &\mapsto \pi^*(x) \smile u_R
\end{align*} 
is an isomorphism for all $j \geq 0$.  Note that, in the orientable case, the Thom class $u_R$ is the image of the integral Thom class $u_\Z$ under the homomorphism $H^{k+1}(E, \partial E; \Z) \to H^{k+1}(E, \partial E; R)$ induced by the homomorphism $\Z \to R$ taking $1 \in \Z$ to the identity of $R$.  In particular, the $\Z_2$-Thom class is the mod-$2$ reduction of the integral Thom class.

The $R$-Euler class $e_R \in H^{k+1}(B; R)$ of the bundles $D^{k+1} \to E \xrightarrow{\pi} B$ and $\sph^k \to \partial E \xrightarrow{\pi} B$ is defined by
$$
e_R = \Phi^{-1}(u_R \smile u_R)
$$
and gives rise to the long exact Gysin sequence
$$
\cdots \to H^{j+k}(\partial E) \to H^j(B) \xrightarrow{\smile e_R} H^{j+ k + 1}(B) \to H^{j+ k + 1}(\partial E) \to \cdots
$$
with coefficients in  $R$.  If $k + 1$ is odd, observe that $2 (u_R \smile u_R) = 0$ and, hence, the Euler class satisfies $2 \, e_R = 0$.  Equivalently, the Euler class may be defined via $e_R = (\pi^*)^{-1} \circ j^*(u_R)$, where $j : E \to (E, \partial E)$ is the (relative) inclusion map.

From now on, the subscript $R$ will be suppressed from the Thom and Euler classes.  In the cases below where $\Z_2$ coefficients need to be highlighted, the notation  $\bar u$ and $\bar e$ will denote the $\Z_2$ classes, respectively.

Before discussing Stiefel-Whitney classes, recall first that the Steenrod squares (see, for example, \cite[Section 4.L]{Hatcher} and \cite[Section 8.1]{MS}) are additive cohomology operations $\Sq^i:H^j(X, A; \Z_2) \to H^{j + i}(X, A; \Z_2)$ which satisfy the following properties for any pair $(X,A)$:

\begin{enumerate}[label=(\roman*), leftmargin=10mm]
\item  $\Sq^0 : H^j(X, A; \Z_2) \to H^j(X, A; \Z_2)$ is the identity homomorphism.

\vspace{1mm}
\item  $\Sq^1 \colon H^j(X, A; \Z_2) \!\to\! H^{j + 1}(X, A; \Z_2)$ is the Bockstein homomorphism associated to the short exact sequence $0\to  \Z_2 \to \Z_4 \to  \Z_2 \to 0$.

\vspace{1mm}
\item  $\Sq^i(x) = x \smile x$ if $\deg(x) = i$.

\vspace{1mm}
\item  $\Sq^i(x) = 0$ if $\deg(x) < i$. 

\vspace{1mm}
\item  (Naturality)  If $f:(X, A)\rightarrow (Y, B)$ is continuous, then 
$$
f^\ast \circ \Sq^i = \Sq^i \circ f^\ast\,.
$$

\vspace{1mm}
\item  (Cartan formula) If $x \smile y$ is defined, then 
$$
\Sq^i(x \smile y) = \sum_{j=0}^i \Sq^j(x) \smile \Sq^{i-j}(y) \,.
$$

\vspace{1mm}
\item  (Adem relations \cite{Adem})  If $i < 2j$, then 
$$
\Sq^i \circ \Sq^j = \sum_{l=0}^{\lfloor i/2\rfloor} \binom{j-l-1}{i-2l} \Sq^{i+j-l} \circ \Sq^l \,.
$$

\end{enumerate}

With the Steenrod squares in hand, the Stiefel-Whitney classes $w_i(\pi) \in H^i(B; \Z_2)$ of the bundles $D^{k+1} \to E \xrightarrow{\pi} B$ and $\sph^k \to \partial E \xrightarrow{\pi} B$ are now defined (see \cite{Hsiang}, \cite{Thom}) via
$$
w_i(\pi) = \Phi^{-1}(\Sq^i(\bar u)) \,.
$$
In particular, observe that $w_0(\pi) = 1$ and that the Euler class satisfies $\bar e = w_{k+1}(\pi)$.  

As in the linear case, the Stiefel-Whitney classes for smooth bundles satisfy the Wu formula \cite{Hsiang}
\begin{equation}
\label{E:Wu}
\Sq^i(w_j(\pi)) = \sum_{l=0}^i \binom{j - i + l - 1}{l} w_{i-l}(\pi) \smile w_{j + l}(\pi) \mod \, 2\,,
\end{equation}
where the binomial coefficient $\binom{a}{b}$ is defined as usual for $a \geq b \geq 0$, to be equal to $1$ if $(a,b) = (-1,0)$, and to equal $0$ otherwise.  

The Wu formula can be used to establish the well-known property that the first non-trivial Stiefel-Whitney class must occur in degree $2^r$, for some $r \geq 0$.  Thus, it is clear that the first non-trivial Stiefel-Whitney class must occur in even degree whenever $w_1(\pi) = 0$.  Furthermore, if both $w_1(\pi) = 0$ and $w_2(\pi) = 0$, then the first non-trivial Stiefel-Whitney class must occur in degree $d = 2^{2+r} \equiv 0$ mod $4$, for some $r \geq 0$, and, by considering $i=2$ and $j = d-1$ in the Wu formula, it follows in this case that $w_{d+1}(\pi) = 0$.

\begin{rem}
It is also possible to define rational Pontryagin classes for general smooth disk and sphere bundles.  Indeed, there is a natural notion of rational Pontryagin class for fiber bundles with fiber $\R^{k+1}$ and structure group  $\Homeo(\R^{k+1})$ (see the introduction of \cite{Weiss}).  This yields rational Pontryagin classes for smooth sphere bundles by considering those $\R^{k+1}$-bundles whose structure group reduces to $\Diff(\sph^k)$ via the homomorphism $\Diff(\sph^k) \to \Homeo(\R^{k+1})$ given by extending radially.  Similarly, by composing with the restriction homomorphism $\Diff(D^{k+1}) \to \Diff(\sph^k)$, that is, by considering those $\R^{k+1}$-bundles whose structure group further reduces to $\Diff(D^{k+1})$, one obtains rational Pontryagin classes for smooth disk bundles.  Note that these rational Pontryagin classes are not automatically trivial, since the inclusion homomorphism $\Or(k+1) \to \Homeo(\R^{k+1})$ (which yields the linear rational Pontryagin classes) factors through the groups $\Diff(D^{k+1})$ and $\Diff(\sph^k)$.  However, in the non-linear setting, Weiss \cite{Weiss} has demonstrated that the rational Pontryagin classes cannot be expected to satisfy many of the properties satisfied by Pontryagin classes of linear bundles.
\end{rem}


\section{Even-dimensional rational homology spheres}
\label{S:evendim}

The goal of this section is to prove Theorem \ref{thm:main}.  Throughout this section, let $M^n$ be an even-dimensional, smooth, closed, simply connected manifold with the rational homology of $\sph^n$ (that is, $M^n$ is an even-dimensional rational homology sphere), such that $M^n$ admits a double disk-bundle decomposition $DB_- \cup_L DB_+$.  Recall that, by \cite[Proposition 4.1]{DGGK}, it may always be assumed that the singular leaves $B_\pm$ are connected and that the regular leaf $L$ is the total space of sphere bundles $\sph^{\ell_\pm}\rightarrow L\rightarrow B_\pm$ over the singular leaves.  By Remark \ref{R:codim}, we have that $\ell_\pm \geq 1$ whenever $B_\pm$ are connected.  

The strategy of the proof is to reduce the problem to two cases: namely, the easy case, where the singular leaves $B_\pm$ have $\dim(B_\pm) = 0$, and the more difficult case, where the regular leaf $L$ has $\dim(L) = \dim(B_-) + \dim(B_+)$.  In the more difficult case, it will be shown that the singular leaves must both be homotopy spheres, from which it may be concluded that the even-dimensional rational homology sphere $M^n$ is actually a homotopy sphere.

\begin{lem}
\label{L:dim_constraints}
Suppose that $M^n$ is a simply connected rational homology sphere of even dimension $n \geq 6$.  Suppose, furthermore, that $M^n$ admits a double disk-bundle decomposition $DB_- \cup_L DB_+$ with $B_\pm$ connected.  Then the bundles $\sph^{\ell_\pm}\rightarrow L\rightarrow B_\pm$ are both orientable and either
\begin{enumerate}
\item $n = \ell_\pm + 1$ and $DB_\pm \cong D^n$, or

\vspace{1mm}
\item 
\label{i:dimsum}
$n = \ell_- + \ell_+ + 1$, with $\dim(B_+) = \ell_- \neq \ell_+ = \dim(B_-)$ and exactly one of $\ell_\pm$ is odd.
\end{enumerate}
\end{lem}

\begin{proof}
By Corollary \ref{C:dims} and since $n$ is even, the homotopy fiber $F$ of the inclusion $L \to M$ has a loop space factor $\Omega \sph^k$ with either $k = n$ or $k = 2n - 1 \equiv 3$ mod $4$.  Moreover, by Corollary \ref{C:torF} and since $n \geq 6$, the homotopy fiber $F$ has torsion-free integral cohomology.  Therefore, it follows from Table \ref{table:cohomF} that the bundles $\sph^{\ell_\pm} \to L \to B_\pm$ are both orientable.  From Table \ref{table:Qtype}, it now follows that $k$ and $\ell_\pm$ must satisfy one of the following:
\begin{enumerate}[label=(\roman*)]
\item \label{i:constraint1}
$k = 3$ and $\ell_\pm = 1$;

\vspace{1mm}
\item \label{i:constraint2}
$k = \ell_- + \ell_+ + 1$ and $1 \leq \ell_\pm$;

\vspace{1mm}
\item \label{i:constraint3}
$k = \ell_\pm + 1$ and $1 < \ell_- = \ell_+$;

\vspace{1mm}
\item \label{i:constraint4}
$k = 7$ and $\ell_\pm = 2$.
\end{enumerate}
In cases \ref{i:constraint1} and \ref{i:constraint4}, the only possibility is $k = 2n-1$, since $n$ is even.  However, in both cases, this forces $n \leq 4$, contradicting the assumption that $n \geq 6$.  Therefore, these cases cannot occur.

In case \ref{i:constraint3}, suppose first that $2n-1 = k = \ell_\pm + 1$.  This implies that $2n-1 = \ell_\pm + 1 \leq \dim(L) + 1 = n$, which is impossible, since $n \geq 6$.  Hence, the only possibility is that $n = k = \ell_\pm + 1$.  This implies that $\dim(B_\pm) = \dim(L) - \ell_\pm = 0$, so the (connected) singular leaves $B_\pm$ must each consist of a single point and, consequently, it may be concluded that $DB_\pm \cong D^n$.  

Finally, in case \ref{i:constraint2}, suppose first that $2n-1 = k = \ell_- + \ell_+ + 1$.  Since $\ell_\pm \leq \dim(L) = n-1$, this is possible only if $\ell_\pm = n - 1$, leading to $\dim(B_\pm) = 0$ and $DB_\pm \cong D^n$, as before.  Thus, suppose instead that $n = k = \ell_- + \ell_+ + 1$.  As $n$ is even, it follows that $n - 1 = \ell_- + \ell_+$ is odd, meaning that $\ell_- \neq \ell_+$, with exactly one of them being odd.  Moreover, the singular leaves have dimensions $\dim(B_\pm) = \dim(L) - \ell_\pm = (n-1) - \ell_\pm = \ell_\mp$ respectively, as desired.
\end{proof}

Case \eqref{i:dimsum} of Lemma \ref{L:dim_constraints} is precisely what happens for the standard double disk-bundle decomposition
$$
\sph^n = (\sph^p \times D^{q+1}) \cup_{\sph^p \x \sph^q} (D^{p+1}\times \sph^q)
$$ 
of an even-dimensional sphere, where $p + q + 1 = n$ and $p\neq q$.  The content of Corollary \ref{cor:homeotype}, established in Corollary \ref{C:corB} below, is that this standard decomposition is, up to homeomorphism,  the model for every non-trivial double disk-bundle decomposition of an even-dimensional sphere of dimension $\geq 6$.

\begin{prop}
\label{P:one_sphere}
Let $M^n$ be a simply connected rational homology sphere of even dimension $n \geq 6$ which admits a double disk-bundle decomposition $DB_- \cup_L DB_+$ with $B_\pm$ connected.  Suppose that $n = \ell_- + \ell_+ +1$ and, without loss of generality, that $\ell_- < \ell_+$.  Then 
\begin{enumerate}
\item $B_+$ is a homotopy sphere, 

\vspace*{1mm}
\item $B_-$ is a simply connected rational homology sphere and

\vspace*{1mm}
\item $L$ has the same integral cohomology ring as the product $\sph^{\ell_- } \x \sph^{\ell_+}$.
\end{enumerate}
\end{prop}

\begin{proof}
By Lemma \ref{L:dim_constraints}, the bundles $\sph^{\ell_\pm} \to L \xrightarrow{\pi_\pm} B_\pm$ are both orientable and the singular leaves have dimensions $\dim(B_\pm) = \ell_\mp$.  Moreover, since the regular leaf $L$ is orientable, being a hypersurface in a simply connected manifold, it follows that $B_\pm$ are also orientable.

Since $n = \ell_- + \ell_+ + 1$ and $\dim(B_\pm) = \ell_\mp$, the isomorphisms \eqref{E:isoms} imply that $H^{j+1}(M, B_\pm) \cong H_{n - j - 1}(B_\mp) = 0$ for all $j \in \{1, \dots, \ell_\mp - 1\}$. Therefore, the long exact sequences 
$$
\cdots \to H^j(M) \to H^j(B_\pm) \to H^{j+1}(M, B_\pm) \to H^{j+1}(M) \to \cdots
$$
for the pairs $(M, B_\pm)$ (contained in the braid diagram \eqref{E:braid}) yield that $H^j(M) \to H^j(B_\pm)$ is surjective and, hence, that $H^j(B_\pm)$ is finite, for all $j \in \{1, \dots, \ell_\mp - 1\}$, since $M^n$ is a rational homology sphere.  In other words, the singular leaves $B_\pm$ are both rational homology spheres.  

Now, since $1 \leq \ell_- < \ell_+$, Proposition \ref{prop:pi1} implies that the homomorphism $\pi_1(\sph^{\ell_-}) \to \pi_1(L)$ in the long exact homotopy sequence for the bundle $\sph^{\ell_-} \to L \to B_-$ must be a surjection and, hence, that $\pi_1(B_-) = 0$.  Moreover, observe that $B_+$ is diffeomorphic to $\sph^1$ whenever $\dim(B_+) = \ell_- = 1$.  On the other hand, if $2 \leq \ell_- < \ell_+$, then Proposition \ref{prop:pi1}, together with the long exact homotopy sequence for the bundle $\sph^{\ell_+} \to L \to B_+$, implies that $L$ and $B_+$ are both simply connected.  Thus, by the classification of surfaces, the resolution of the Poincar\'e Conjecture and the classification of $4$-manifolds \cite{Fr}, the simply connected rational homology sphere $B_+$ is a homotopy sphere whenever $2 \leq \ell_- \leq 4$.  Hence, it may be assumed from now on, without loss of generality, that $5 \leq \ell_- < \ell_+$.

Let $e_+\in H^{\ell_+ + 1}(B_+)$ denote the Euler classes of the bundle $\sph^{\ell_+} \to L \xrightarrow{\pi_+} B_+$.  Since $\dim(B_+) = \ell_- < \ell_+ + 1$ by hypothesis, it follows that $e_+ = 0$.  Therefore, the Gysin sequence for $\sph^{\ell_+} \to L \to B_+$ decomposes into short exact sequences
$$
0 \to H^j(B_+) \xrightarrow{\pi_+^*} H^j(L) \to H^{j - \ell_+}(B_+) \to 0 \,.
$$
Again using the fact that $\dim(B_+) = \ell_- < \ell_+$, it now follows that $H^j(B_+)$ and $H^{j - \ell_+}(B_+)$ cannot be non-trivial simultaneously and, therefore, that
$$
H^j(L) \cong 
\begin{cases}
H^j(B_+), & \text{ if } j \in \{0, \dots, \ell_-\}, \\
H^{j - \ell_+}(B_+), & \text{ if } j \in \{\ell_+, \dots, n - 1\}, \\
0, & \text{ otherwise.}
\end{cases}
$$
This implies, in particular,  that $H^j(L) \cong \Z$ for all $j \in \{0, \ell_\pm, n-1\}$, since $B_+$ is connected and orientable.  In fact, more can be said.  Since $H^j(DB_+) \cong H^j(B_+) = 0$ for all $j > \ell_-$, the (relative) fiber inclusion $\iota : (D^{\ell_+ + 1}, \sph^{\ell_+}) \to (DB_+, L)$ induces, by naturality of the long exact sequence of a pair and the existence of a Thom class $u_+ \in H^{\ell_+ + 1}(DB_+, L)$, a commutative diagram
$$
\xymatrix{
H^{\ell_+ }(\sph^{\ell_+}) \ar[r]_(0.4)\cong & H^{\ell_+ + 1}(D^{\ell_+ + 1}, \sph^{\ell_+}) \\
H^{\ell_+ }(L) \ar[r]^(0.4){\cong} \ar[u]^{\iota^*}_\cong & H^{\ell_+ + 1}(DB_+, L) \ar[u]_{\iota^*}^\cong
}
$$
of isomorphisms between free abelian groups of rank one.  Therefore, a generator of $H^{\ell_+}(L)$ restricts to a generator of $H^{\ell_+}(\sph^{\ell_+})$ for each fiber of the bundle $\sph^{\ell_+} \to L \to B_+$.  By the Leray-Hirsch Theorem, it now follows that $H^*(L)$ is a free $H^*(B_+)$-module with basis given by a generator of $H^0(L) \cong H^0(\sph^{\ell_+})$ and a generator of $H^{\ell_+}(L) \cong H^{\ell_+}(\sph^{\ell_+})$, where the scalar multiplication in $H^*(L)$ is given by 
\begin{align*}
H^*(B_+) \x H^*(L) &\to H^*(L) \\
(\beta, x) &\mapsto \pi_+^*(\beta) \smile x \,.
\end{align*}
In particular, since $\pi_+^* : H^j(B_+) \to H^j(L)$ is an isomorphism for all $j \in \{0, \dots, \ell_-\}$ and since $H^{\ell_+}(L) \cong \Z$, it follows that, for all $j \in \{0, \dots, \ell_-\}$, the cup product $x \smile y \in H^{j + \ell_+}(L)$ must be non-trivial whenever $x \in H^j(L)$ and $y \in H^{\ell_+}(L)$ are both non-trivial.
 
 Suppose now that the simply connected rational homology sphere $B_+$ is not a homotopy sphere, with $5 \leq \dim(B_+) = \ell_-$.  Then, by Poincar\'e duality and the universal coefficient theorem, there exists some $j_0 \in \{3, \dots, \ell_- - 2\}$ such that $H^j(L) \cong H^j(B_+) = 0$ for all $j \in \{1, \dots, j_0 - 1\}$ and such that $H^{j_0}(L) \cong H^{j_0}(B_+)$ is a non-trivial finite group.
 
On the other hand, since the second singular leaf $B_-$ is a simply connected rational homology sphere of dimension $\dim(B_-) = \ell_+ > \ell_-$, it follows from the Gysin sequence for the bundle $\sph^{\ell_-} \to L \xrightarrow{\pi_-} B_-$ that $\pi_-^* : H^{j_0}(B_-) \to H^{j_0}(L)$ is an isomorphism and that $\pi_-^* : H^{\ell_+}(B_-) \to H^{\ell_+}(L)$ is injective.  Let $x_- \in H^{j_0}(B_-) \cong H^{j_0}(L) \cong H^{j_0}(B_+)$ and $y_- \in H^{\ell_+}(B_-) \cong \Z$ be non-trivial cohomology classes, and let $x = \pi_-^*(x_-) \in H^{j_0}(L)$ and $y = \pi_-^*(y_-) \in H^{\ell_+}(L)$ be the corresponding non-trivial cohomology classes in $L$.  From our observations above, it follows that $x \smile y \in H^{j_0 + \ell_+}(L)$ must be non-trivial.  However, since $\dim(B_-) = \ell_+$, we have $x_- \smile y_- = 0$ and, hence, that
$$
0 = \pi_-^*(x_- \smile y_-) = \pi_-^*(x_-) \smile \pi_-^*(y_-) = x \smile y\,.
$$
This contradiction implies that $B_+$ must be a homotopy sphere, as claimed.

Finally, the fact that $H^*(L)$ has the same ring structure as $H^*(\sph^{\ell_- } \x \sph^{\ell_+})$ now follows from Poincar\'e duality \cite[Corollary 3.39]{Hatcher}, since the non-trivial cohomology groups of $L$ occur in degrees $0, \ell_\pm, n-1$ and are free abelian of rank one, and from the fact that the only scenario in which the cup product $\smile : H^{\ell_-}(L) \x H^{\ell_-}(L) \to H^{2 \ell_-}(L)$ could possibly be non-trivial is when $\ell_+ = 2 \ell_-$, which would imply that $\ell_-$ is odd and, hence, that the square of a generator of $H^{\ell_-}(L) \cong \Z$ is a torsion element in $H^{2 \ell_-}(L) = H^{\ell_+}(L) \cong \Z$, which is clearly impossible.
\end{proof}

In the context of Proposition \ref{P:one_sphere}, it turns out that the simply connected rational homology sphere $B_-$ is also a homotopy sphere.  However, in order to prove this, it is first necessary to establish a preliminary lemma.  Observe that the sphere bundle
$$
\sph^2 \cong \SO(3)/\SO(2) \to \SU(3)/\SO(2) \to \SU(3)/\SO(3)
$$ 
that arises from the inclusions $\SO(2) \In \SO(3) \In \SU(3)$ gives an example of a sphere bundle over the Wu manifold $\SU(3)/\SO(3)$ (a rational homology $5$-sphere) where the total space has the cohomology ring of a product of two spheres.  That the total space $\SU(3)/\SO(2)$ has such a cohomology ring follows from the observation that it also occurs as the total space of the sphere bundle $\sph^2 \cong \SU(2)/\SO(2) \to \SU(3)/\SO(2) \to \SU(3)/\SU(2) \cong \sph^5$ with trivial Euler class coming from the inclusions $\SO(2) \In \SU(2) \In \SU(3)$.  The point of the following lemma (to be applied in the proof of Theorem \ref{T:two_spheres}) is that such a bundle over the Wu manifold $\SU(3)/\SO(3)$ is inadmissable when considering double disk-bundle decompositions of even-dimensional rational homology spheres.  

\begin{lem}
\label{L:noWu}
Let $\sph^5 \to L \to \sph^2$ be a smooth sphere bundle.  Then $L$ cannot be the total space of a smooth sphere bundle $\sph^2 \to L \to B$ over the Wu manifold $B = \SU(3)/\SO(3)$.
\end{lem}

\begin{proof}
Observe first that, by a similar argument to that in the proof of Proposition \ref{P:one_sphere}, the Euler class of the bundle $\sph^5 \to L \to \sph^2$ is trivial and, hence, the total space $L$ must have the same cohomology ring as the product $\sph^2 \x \sph^5$.  Moreover, from the long exact homotopy sequence for the bundle $\sph^5 \to L \to \sph^2$, it follows that $\pi_j(L) \cong \pi_j(\sph^2)$ for $j \leq 4$, that there is a short exact sequence $0 \to \Z \cong \pi_5(\sph^5) \to \pi_5(L) \to \pi_5(\sph^2) \cong \Z_2 \to 0$, and that $\pi_6(L)$ surjects onto $\pi_6(\sph^2) \cong \Z_{12}$.

Suppose now that $L$ is also the total space of a smooth sphere bundle $\sph^2 \to L \to B$ over the Wu manifold $B = \SU(3)/\SO(3)$.  The Hurewicz theorem yields $\pi_2(B) \cong H_2(B) \cong \Z_2$, while Propositions 4.1 and 4.2 of \cite{Ka} (see also \cite[Theorem 5.6]{PR}) yield $\pi_3(B) \cong \Z_4$, $\pi_4(B) \cong 0$, $\pi_5(B) \cong \Z \oplus \Z_2$ and $\pi_6(B) \cong \Z_2$.  Therefore, from the long exact homotopy sequence for the bundle $\sph^2 \to L \to B$ it follows that the homomorphism $\pi_5(L) \to \pi_5(B) \cong \Z \oplus \Z_2$ is surjective.  Since it was shown above that $\pi_5(L)$ fits into a short exact sequence $0 \to \Z \to \pi_5(L) \to \Z_2 \to 0$, it may be concluded that $\pi_5(L) \to \pi_5(B) \cong \Z \oplus \Z_2$ is, in fact, an isomorphism.  Hence, in the long exact sequence, the homomorphism $\Z_2 \cong \pi_6(B) \to \pi_5(\sph^2) \cong \Z_2$ must also be an isomorphism, which, in turn, implies that the homomorphism $\Z_{12} \cong \pi_6(\sph^2) \to \pi_6(L)$ is surjective.  Therefore, together with the previous paragraph, it now follows that $\pi_6(L) \cong \pi_6(\sph^2) \cong \Z_{12}$.

The remainder of the proof, modelled on an argument in Section 1.2 of \cite{Kr}, will establish that the isomorphism $\pi_6(L) \cong \Z_{12}$ is impossible.  Recall that $L$ is simply connected, with $H^*(L) \cong H^*(\sph^2 \x \sph^5)$.  Let $\sph^1 \to X \to L$ be the principal circle bundle over $L$ associated to a generator of $H^2(L) \cong \Z$, so that the total space $X$ is a simply connected $8$-manifold.  By the Gysin sequence, together with Poincar\'e duality \cite[Corollary 3.39]{Hatcher}, it follows that $H^*(X) \cong H^*(\sph^3 \x \sph^5)$, while the long exact homotopy sequence for the bundle $\sph^1 \to X \to L$ yields $\pi_2(X) = 0$ and $\pi_j (X) \cong \pi_j(L)$, for all $j \geq 3$. 

Now, by Proposition 4C.1 of \cite{Hatcher}, there exists a CW-complex $Y$ which is homotopy equivalent to $X$, such that
$$
Y = \sph^3 \cup_{\alpha_4} e_5 \cup_{\alpha_7} e_8 \,,
$$
where $\alpha_k : \sph^{k} \cong \partial e_{k+1} \to Y^{(k)}$ denotes the map attaching the $(k+1)$-cell $e_{k+1}$, $k \in \{4,7\}$, to the $k$-skeleton $Y^{(k)}$ of $Y$.  Moreover, recall that the pair $(Y, Y^{(k)})$ is $k$-connected, by Corollary 4.12 of \cite{Hatcher}, meaning that the inclusion $\iota : Y^{(k)} \to Y$ induces isomorphisms $\iota_* : \pi_j(Y^{(k)}) \to \pi_j(Y)$, for all $j \leq k-1$, and a surjective homomorphism $\iota_* : \pi_{k}(Y^{(k)}) \to \pi_{k}(Y)$.

Since $Y^{(4)} = \sph^3$, $\pi_4(\sph^3) \cong \Z_2$ and $\pi_4(Y) \cong \pi_4(X) \cong \pi_4(L) \cong \Z_2$, it follows that the surjective homomorphism $\iota_* : \pi_{4}(Y^{(4)}) \to \pi_{4}(Y)$ is actually an isomorphism.  This implies that the attaching map $\alpha_4$ must be null-homotopic, as its homotopy class $[\alpha_4] \in \pi_{4}(Y^{(4)})$ would otherwise be mapped under $\iota_*$ to the non-trivial class in $\pi_4(Y)$, contradicting the fact that $\iota \circ \alpha_4$ clearly extends to a map of the $5$-cell $e_5$ into Y.  Hence, since $\alpha_4$ is null-homotopic, it may be assumed without loss of generality that $Y = (\sph^3 \vee \sph^5) \cup_{\alpha_7} e_8$.

Finally, since $(Y, Y^{(7)})$ is $7$-connected, where $Y^{(7)} = \sph^3 \vee \sph^5$, and since the pair $(\sph^3 \x \sph^5, \sph^3 \vee \sph^5)$ is also $7$-connected (by Exercise 17 in Chapter 4.1 of \cite{Hatcher}), it may be concluded that
$$
\pi_6(L) \cong \pi_6(X) \cong \pi_6(Y) \cong \pi_6(Y^{(7)}) \cong \pi_6(\sph^3 \x \sph^5) \cong \Z_{12} \oplus \Z_2 \,,
$$
contradicting the previous assertion that $\pi_6(L) \cong \Z_{12}$.
\end{proof}

\begin{thm}
\label{T:two_spheres}
Let $M^n$ be a simply connected rational homology sphere of even dimension $n \geq 6$ which admits a double disk-bundle decomposition $DB_- \cup_L DB_+$ with $B_\pm$ connected.  Suppose that $n = \ell_- + \ell_+ +1$ and, without loss of generality, that $\ell_- < \ell_+$.  Then $B_-$ is a homotopy sphere.
\end{thm}

\begin{proof}
Recall from Lemma \ref{L:dim_constraints} that the bundles $\sph^{\ell_\pm} \to L \xrightarrow{\pi_\pm} B_\pm$ are both orientable, and from Proposition \ref{P:one_sphere} that $B_+$ is a homotopy sphere, that $B_-$ is simply connected rational homology sphere and that $L$ has the same cohomology ring as the product $\sph^{\ell_-} \x \sph^{\ell_+}$.

It is sufficient to show that the Euler class $e_- \in H^{\ell_- + 1}(B_-)$ of the bundle $\sph^{\ell_-} \to L \xrightarrow{\pi_-} B_-$ is trivial.  Indeed, if $e_- = 0$, then the Gysin sequence for the bundle decomposes into short exact sequences
$$
0 \to H^j(B_-) \xrightarrow{\pi_-^*} H^j(L) \to H^{j - \ell_-}(B_-) \to 0
$$
and, since $H^j(L)$ is free abelian for all $j$, the injective homomorphisms $H^j(B_-) \to H^j(L)$ yield $H^j(B_-) = 0$ for all $j \in \{1, \dots, \ell_+ - 1\}$, as desired.

Assume, therefore, that the Euler class $e_- \in H^{\ell_- + 1}(B_-)$ is non-trivial.  Observe first that, for all $j \in \{1, \dots, \ell_-\}$, the Gysin sequence yields injective homomorphisms $\pi_-^* : H^j(B_-) \to H^{j}(L)$, from which it may be concluded that $H^j(B_-) = 0$ for all $j \in \{1, \dots, \ell_-\}$, since $H^{j}(L)$ is free abelian.  Moreover, it also follows that $\ell_+ = \dim(B_-) > \ell_- + 1$, since, if $\ell_+ = \ell_- + 1$ were true, the exactness of the Gysin sequence and the fact that $L$ has the same cohomology ring as the product $\sph^{\ell_-} \x \sph^{\ell_+}$ would together lead to the homomorphism $H^0(B_-) \cong \Z \to H^{\ell_- + 1}(B_-) = H^{\ell_+}(B_-) \cong \Z$ given by cupping with $e_-$  being trivial, contradicting the fact that $e_-$ is both non-trivial and in the image of this homomorphism.

The inequality $\ell_+ > \ell_- + 1$ further implies that $\ell_+ = \dim(B_-) > 2 \ell_-$, since, if $\ell_+ \leq 2 \ell_-$ were true, Poincar\'e duality and the universal coefficient theorem would yield the isomorphisms $H^{\ell_- + 1}(B_-) \cong H_{\ell_+ - \ell_- - 1}(B_-) \cong H^{\ell_+ - \ell_-}(B_-)$ of finite groups, whereas the inequalities $1 < \ell_+ - \ell_- \leq \ell_-$ imply $H^{\ell_+ - \ell_-}(B_-) = 0$, thus again contradicting the assumption that $e_- \in H^{\ell_- + 1}(B_-)$ is non-trivial.

Therefore, without loss of generality, it may be assumed from now on that $\dim(B_-) = \ell_+ \geq 2 \ell_- + 1$.  In other words, there exist some $m \in \N$, $m \geq 2$, and some $r \in \{-1, 0, 1, \dots, \ell_- - 1\}$, such that $\ell_+ = m(\ell_- + 1) + r$.

Suppose first that $r \in \{0, \dots, \ell_- - 1\}$.  Then, since $\dim(B_-) = \ell_+ = m(\ell_- + 1) + r$, it follows  that $H^{m(\ell_- + 1)}(B_-) \cong \Z$, whenever $r = 0$, and, via Poincar\'e duality and the universal coefficient theorem, that $H^{m(\ell_- + 1)}(B_-) \cong H_r(B_-) \cong H^{r+1}(B_-) = 0$, whenever $r \in \{1, \dots, \ell_- - 1\}$.

Recall that the inequalities $0 < (m-1)(\ell_- + 1) < \ell_+ = \dim(B_-)$ imply that $H^{(m-1)(\ell_- + 1)}(B_-)$ is finite, while the inequalities $\ell_- < m(\ell_- + 1) - 1 < m(\ell_- + 1) + r = \ell_+$ imply that $H^{m(\ell_- + 1) - 1}(L) = 0$.  By combining these observations with the fact that $H^{m(\ell_- + 1)}(B_-)$ is free abelian of rank $\leq 1$, it follows from the portion 
$$
\cdots \to H^{m(\ell_- + 1) - 1}(L) \to H^{(m-1)(\ell_- + 1)}(B_-) \to H^{m(\ell_- + 1)}(B_-) \to \cdots
$$
of the Gysin sequence that $H^{(m-1)(\ell_- + 1)}(B_-) = 0$.  If $m = 2$, then this already contradicts the assumption that $e_- \neq 0$.   On the other hand, if $m \geq 3$, then induction (via an analogous argument for the induction step) shows that $H^{j(\ell_- + 1)}(B_-) = 0$ for all $j \in \{1, \dots, m-1\}$, again yielding the contradiction $H^{\ell_- + 1}(B_-) = 0$.

It remains, therefore, to examine the case $r = -1$: that is, to show that assuming that $e_- \in H^{\ell_- + 1}(B_-)$ is non-trivial leads to a contradiction whenever $\ell_+ = m(\ell_- + 1) -1$, for some $m \in \N$, $m \geq 2$.  First observe that $\dim(L) = \ell_- + \ell_+ = (m-1)(\ell_- + 1) + 2 \ell_-$ being odd implies that $\ell_-$ and $m$ must both be even and, hence, that $\ell_+ = m(\ell_- + 1) - 1$ is odd.  Now, since $H^{\ell_-}(B_-) = 0$ and $H^{\ell_- + 1}(L) = 0$ (as $\ell_- < \ell_- + 1 < \ell_+$), the Gysin sequence yields a short exact sequence
$$
0 \to H^{\ell_-}(L) \cong \Z \to H^0 (B_-) \cong \Z \xrightarrow{\smile e_-} H^{\ell_- + 1}(B_-) \to 0 \,.
$$
Moreover, as discussed in Section \ref{S:nonlin}, the fiber of the bundle $\sph^{\ell_-} \to L \to B_-$ being even-dimensional implies that the Euler class $e_- \in H^{\ell_- + 1}(B_-)$ must be torsion of order two.  Since $e_- \neq 0$ by assumption, it follows from the short exact sequence above that $H^{\ell_- + 1}(B_-) = \<e_-\>/\<2 e_-\> \cong \Z_2$.

Since $H^*(L) \cong H^*(\sph^{\ell_-} \x \sph^{\ell_+})$, it now follows from the Gysin sequence for $\sph^{\ell_-} \to L \to B_-$ that the cup product $\smile e_- : H^j(B_-) \to H^{j + \ell_- + 1}(B_-)$ is an isomorphism for all $j \in \{1, \dots, \ell_+ - \ell_- - 2\}$.  Given that $\dim(B_-) = \ell_+ = m(\ell_- + 1) - 1$, that $H^j(B_-) = 0$, for all $j \in \{1, \dots, \ell_-\}$, and that $H^{\ell_- + 1}(B_-) = \<e_-\>/\<2 e_-\> \cong \Z_2$, it thus follows that the cohomology ring of $B_-$ is given by
$$
H^*(B_-) \cong \<e_-, x\> / \< 2 e_-, e_-^m, x^2, x e_-\> \,,
$$
where $\deg(x) = \ell_+$.  In particular, if $0 < j < \dim(B_-)$, then $H^j(B_-)$ is non-trivial if and only if $j \equiv 0$ mod $\ell_- + 1$, in which case $H^j(B_-) \cong \Z_2$ is generated by $e_-^j$.

Via the universal coefficient theorem, it follows that the $\Z_2$-cohomology groups of $B_-$ are given by 
$$
H^j(B_-; \Z_2) = \begin{cases}
0, & \text{ if } j \in \{1, \dots, \ell_- - 1\} \mod \ell_- + 1, \\
\Z_2, & \text{ if } j \in \{\ell_-, \ell_- + 1\} \mod \ell_- + 1.
\end{cases}
$$
If $w_j^- = w_j(\pi_-) \in H^j(B_-; \Z_2)$ denote the Stiefel-Whitney classes of the bundle $\sph^{\ell_-} \to L \xrightarrow{\pi_-} B_-$, then, in particular, the top Stiefel-Whitney class $w_{\ell_- + 1}^- \in H^{\ell_- + 1}(B_-; \Z_2) \cong \Z_2$ is given by the mod-$2$ reduction of the Euler class $e_-$ and, hence, is non-trivial.  Moreover, it follows that the group $H^{k(\ell_- + 1)}(B_-; \Z_2) \cong \Z_2$ is generated by $(w_{\ell_- + 1}^-)^k$, for all $k \in \{1, \dots, m-1\}$.  Now, since $\ell_- + 1 \geq 3$ is odd, the Wu formula \eqref{E:Wu} yields the existence of some $j_0 \in \{1, \dots, \ell_-\}$ of the form $j_0 = 2^d$ such that the Stiefel-Whitney class $w_{j_0}^- \in H^{j_0}(B_-; \Z_2)$ is non-trivial.  Clearly, the only possibility is that $j_0 = \ell_-$ in this situation.  However, this implies that $\ell_- = 2$, since, otherwise, $\ell_-$ would be divisible by $4$ and, consequently, the Wu formula would imply that $w_{\ell_- + 1}^- = 0$ (see the end of Section \ref{S:nonlin}), a contradiction.

Therefore, it may be assumed that $\ell_- = 2$ and $\ell_+ = 3m -1$, for some even $m \geq 2$.  From the Wu formula \eqref{E:Wu}, it follows easily that $\Sq^1(w_2^-) = w_3^-$ and, since $(w_2^-)^2 \in H^4(B_-; \Z_2) = 0$, that
$$
0 = \Sq^2\left( (w_2^-)^2 \right) =  \left(\Sq^1(w_2^-)\right)^2 + 2 \, w_2^- \smile \Sq^2(w_2^-)
= (w_3^-)^2\,.
$$
However, this is a contradiction unless $m=2$, since $(w_3^-)^k$ generates the group $H^{3k}(B_-; \Z_2) \cong \Z_2$, for all $k \in \{1, \dots, m-1\}$.  Therefore, the simply connected rational homology sphere $B_-$ has $\dim(B_-) = \ell_+ = 5$ and $H^3(B_-) \cong \Z_2$.  From the Barden-Smale classification of simply connected $5$-manifolds \cite{Ba,Sm3}, the only possibility is that $B_-$ is diffeomorphic to the Wu manifold $\SU(3)/\SO(3)$.

In summary, the only possible case with non-trivial Euler class $e_-$ is the case where $L$ is simultaneously the total space of a sphere bundle $\sph^5 \to L \to \sph^2$ and a sphere bundle $\sph^2 \to L \to B_- \cong \SU(3)/\SO(3)$.  However, by Lemma \ref{L:noWu}, such a scenario is impossible.
\end{proof}

If a rational homology sphere (of any dimension) decomposes as the union of two disk bundles over homotopy spheres, then this has strong topological implications. 

\begin{proposition}
\label{prop:recogNEW} 
Let $M^n$ be a simply connected rational homology sphere which admits a double disk-bundle decomposition $DB_- \cup_L DB_+$, such that the singular leaves $B_\pm$ are connected homotopy spheres of different dimensions.  Then $M^n$ is homeomorphic to a sphere.
\end{proposition}

\begin{proof}
Recall that the bundles $\sph^{\ell_\pm} \to L \xrightarrow{\pi_\pm} B_\pm$ satisfy 
$$
1 \leq \ell_\pm = \dim(L) - \dim(B_\pm) \leq \dim(L) \,, 
$$
since $B_\pm$ are connected.  It is then obvious that $\dim(B_-) \neq \dim(B_+)$ implies $\ell_- \neq \ell_+$.  Moreover, if one of $\ell_\pm$ were equal to $\dim(L) = n-1$, it would follow that the corresponding connected homotopy sphere $B_\pm$ would consist of a single point, which is nonsense.  Therefore, the hypothesis $\dim(B_-) \neq \dim(B_+)$ implies that $n \geq 4$.  To show that $M^n$ is homeomorphic to a sphere, it is sufficient to show that $H^j(M) = 0$ for all $j \in \{1, \dots, n-1\}$ since, in that case, $M^n$ would be a simply connected integral homology sphere.  It would then follow from Whitehead's Theorem that $M^n$ is a homotopy sphere and, hence, by \cite{Fr} and Smale's resolution of the generalised Poincar\'e Conjecture \cite{Sm2}, that $M^n$ is homeomorphic to a sphere.

Since $M^n$ is simply connected, it is trivially true that $H^1(M) = 0$.  Moreover, since $M^n$ is a rational homology sphere and $B_\pm$ are homotopy spheres, it follows that $H^{j}(M)$ is finite and $H^{j}(B_\pm)$ are free abelian for all $j \in \{1, \dots, n-1\}$.  Thus, for $j \in \{1, \dots, n-1\}$, all homomorphisms $H^{j}(M) \to H^{j}(B_\pm)$ in the braid diagram \eqref{E:braid} are trivial.  Hence, for each $j \in \{2, \dots, n-1\}$ there exists a short exact sequence
\beq
\label{E:braidSES}
0 \to H^{j-1}(B_\mp) \to H^{j}(M, B_\mp) \to H^{j}(M) \to 0 \,.
\eeq
If, for each $j \in \{2, \dots, n-1\}$, at least one of $H^{j}(M, B_\mp)$ is trivial, then it is then clear that $H^{j}(M) = 0$.  Assume, therefore, that this is not the case and that there exists some $j_0 \in \{2, \dots, n-1\}$ such that $H^{j_0}(M, B_\mp)$ are both non-trivial.  By the isomorphisms \eqref{E:isoms}, this implies that $H_{n - j_0}(B_\pm) \cong H^{j_0}(M, B_\mp)$ are both non-trivial.  Clearly, the inequalies $2 \leq j_0 \leq n-1$ imply that $1 \leq n - j_0$.  Thus, as $B_\pm$  homotopy spheres, the only possibility is that $n - j_0 = \dim(B_\pm)$, with $H^{j_0}(M, B_\mp) \cong H_{n - j_0}(B_\pm) \cong \Z$,  contradicting the fact that $\dim(B_-) \neq \dim(B_+)$.  
\end{proof}

All the ingredients are now in place to prove Theorem \ref{thm:main}

\begin{thm}
\label{T:evendim} 
Let $M^n$ be an even-dimensional, simply connected rational homology sphere which admits a double disk-bundle decomposition.  Then $M^n$ is homeomorphic to a sphere.
\end{thm}

\begin{proof}
Let $DB_- \cup_L DB_+$ be a double disk-bundle decomposition for $M^n$ and assume, without loss of generality, that the singular leaves $B_\pm$ are connected and, by Remark \ref{R:codim}, that the sphere bundles $\sph^{\ell_\pm} \to L \to B_\pm$ satisfy $\ell_\pm \geq 1$.

If $n=2$, the result follows trivially from the classification of surfaces, while, if $n =4$, it follows from \cite{GeRa}.  Assume, therefore, that $n \geq 6$.

By Lemma \ref{L:dim_constraints}, there are two cases to consider: the case where $\ell_\pm = n-1$ and the case where $n-1 = \ell_- + \ell_+$.  In the first case, it is clear that $B_\pm$ must be points and, hence, that each of the disk bundles $DB_\pm$ is simply a copy of the $n$-disk.  Therefore, in this case $M^n$ is the union of two disks glued along their boundaries and, consequently, homeomorphic to a sphere.

In the case where $n-1 = \ell_- + \ell_+$, it follows from Proposition \ref{P:one_sphere} and Theorem \ref{T:two_spheres} that the singular leaves $B_\pm$ are homotopy spheres.  By Proposition \ref{prop:recogNEW}, this is enough to conclude that $M^n$ is homeomorphic to a sphere, as desired.
\end{proof}

It remains only to prove Corollary \ref{cor:homeotype}.

\begin{corollary}
\label{C:corB}
Suppose $M^{n}$ is an even-dimensional homotopy sphere which admits a double disk bundle decomposition $DB_+\cup_L DB_-$, where $B_\pm$ are both connected.  Then one of the following occurs:

\begin{enumerate}
\item $B_\pm$ are points and $L$ is diffeomorphic to $\sph^{n-1}$.

\item $B_\pm$ are homotopy spheres  and $L$ is homeomorphic to $B_- \times B_+$.

\item $n=4$, $B_\pm$ are diffeomorphic to $\mathbb{R}P^2$ and $L$ is diffeomorphic to $\sph^3/Q_8$.
\end{enumerate}
\end{corollary}

\begin{proof}  Recall that, by Remark \ref{R:codim}, the sphere bundles $\sph^{\ell_\pm} \to L \to B_\pm$ satisfy $\ell_\pm \geq 1$.  Thus, when $n=2$, it is obvious from the identity $n-1 = \dim(L) = \ell_\pm + \dim(B_\pm)$ that case (a) must occur.  The fact that one of (a), (b) or (c) occurs when $n=4$ follows from \cite[Corollary 3.4]{GeRa}.

Assume, therefore, that $n\geq 6$ and suppose that case (a) does not occur.  It follows from Lemma \ref{L:dim_constraints}, Proposition \ref{P:one_sphere} and Theorem \ref{T:two_spheres} that $n = \ell_- + \ell_+ + 1$, that $B_\pm$ are homotopy spheres of dimension $\dim(B_\pm) = \ell_\mp$ respectively and that $L$ has the integral cohomology ring of $\sph^{\ell_-}\times \sph^{\ell_+}$.  It remains to show that $L$ is homeomorphic to $\sph^{\ell_-} \times \sph^{\ell_+}$.  To that end, assume, without loss of generality, that $\ell_- < \ell_+$.

If $\ell_- = 1$ then, since $\ell_+ + \ell_- + 1 = n\geq 6$, it follows that $\dim(B_-) = \ell_+\geq 4$ and, hence, that $B_-$ is $3$-connected.  Since $\Diff(S^1)$ deformation retracts to $\Or(2)$, it may be assumed that the bundle $\sph^1 \to L \to B_-$ is linear.  Now, since $B_-$ is simply connected, the circle bundle $\sph^1 \to L \to B_-$ is automatically orientable and, hence, principal.  As such, it is classified by its Euler class, which must vanish since $H^2(B_-) = 0$.  It follows that the bundle is trivial and, consequently, that $L$ is diffeomorphic to $\sph^1 \times B_-$.

Assume, therefore, that $\dim(B_+) = \ell_- > 1$.  Recall from Section 6.7 and Theorem 12.2 of \cite{Steenrod} that the (smooth) disk bundle $DB_+ \xrightarrow{\pi_+} B_+$ admits a smooth section  $\rho : B_+ \to DB_+$ whose image lies in the interior of $DB_+$.  Then, by the identity $\pi_+ \circ \rho = \id_{B_+}$, the map $\rho$ is an injective immersion.  Since $B_+$ is compact, the image $\rho(B_+)$ is an embedded submanifold of $DB_+$ and, hence, of $M$.  In particular, since $M$ and $B_+$ are homotopy spheres, it follows from Corollary IX.1.4, Theorem IX.7.2 and Corollary IX.8.6 of \cite{Kos} that the normal bundle $\nu B_+$ of $\rho(B_+) \In DB_+ \In M$ is trivial.  Thus, the associated (linear) normal disk bundle $D_\nu B_+ \to \rho(B_+)$ is also trivial, with boundary $\partial D_\nu B_+$ diffeomorphic to $\sph^{\ell_-} \x B_+$.

It now suffices to prove that the total space $D_\nu B_+ \cong D^{\ell_- + 1} \x B_+$ of the normal disk bundle $D_\nu B_+ \to \rho(B_+)$ is diffeomorphic to the total space $D B_+$ of the smooth disk bundle $DB_+ \to B_+$.  To this end, recall that, by Whitehead's Theorem, the bundle projection $\pi_+ : DB_+ \to B_+$ is a homotopy equivalence.  Then it follows from $\pi_+ \circ \rho = \id_{B_+}$ that $\rho_* : H_*(B_+) \to H_*(DB_+)$ is an isomorphism and, hence, from the long exact sequence of the pair $(DB_+,\rho(B_+))$, that $H_\ast(DB_+,\rho(B_+)) = 0$.  Moreover, since $B_+$ is a homotopy sphere of dimension $\dim(B_+) = \ell_- > 1$, it follows from homotopy equivalence that $DB_+$ is simply connected and, from $\ell_+ > \ell_-$ and the long exact homotopy sequence for the bundle $\sph^{\ell_+} \to L \to B_+$, that $\partial D B_+ \cong L$ is also simply connected.  Finally, observe that $\dim(DB_+) = \dim(M) \geq 6$ and $\dim(B_+) + 3 = \ell_- + 3 < \ell_- + \ell_+ + 1 = \dim(DB_+)$.  Thus, the triple $(DB_+, \partial DB_+, \rho(B_+))$ satisfies all of the hypotheses of the Disk Bundle Theorem \cite[Theorem VII.4.4]{Kos} and, hence, it follows that $DB_+$ is diffeomorphic to $D_\nu B_+$, as required.
\end{proof}


\section{Highly connected rational homology spheres}
\label{S:highly_connected}

The goal of this section is to prove Theorem \ref{thm:main2}.  To this end, suppose that $M^n$ is an $(m-1)$-connected, rational homology sphere $M^n$ of dimension $n = 2m+1$, $m \geq 2$.   Thus, it follows that $H^j(M^n) = 0$ for all degrees $j \in \{1, \dots, m \}$ and that $H^{m+1}(M^n) $ is a finite group.  We show that $H^{m+1}(M)$ is a finite cyclic group whenever $M^n$ admits a double disk-bundle decomposition.  

Throughout this section, it will be assumed that $M^n$ admits a double disk-bundle decomposition $DB_- \cup_L DB_+$ with $B_\pm$ connected and $\codim(B_\pm) = \ell_\pm + 1$.  Recall that the closed, smooth, codimension-one submanifold $L = \partial DB_\pm$ is orientable by \cite[p.~107]{Hi} and denote, as usual, the homotopy fiber of the inclusion $L \to M$ by $F$.  

Let $(E_j, d_j)$ denote the spectral sequence for the homotopy fibration $F \to L \to M^n$, where it will be convenient to use the notation $d_j^{p,q}$ to denote the differential $d_j: E_j^{p,q} \to E_j^{p+j, q+1-j}$ on the $E_j$-page.  For each $k \geq 0$, recall that there is a filtration 
$$
0 = A^{k}_{k+1} \In A^k_k \In \cdots \In A^k_0 = H^k(L)
$$
of $H^k(L)$, where $E_\infty^{p, k - p} \cong A^k_p/A^k_{p+1}$.  

At this point, it is convenient to prove a general lemma about spectral sequences using the notation above.
\begin{lem}
\label{L:SES}
If $E_\infty^{i,k-i}$ and $E_\infty^{j, k-j}$, $0 \leq i < j \leq k$, are the only possible non-trivial entries along the diagonal $\{ E_\infty^{p,q} \mid p + q = k \}$, then there exists a short exact sequence 
$$
0 \to E_\infty^{j, k-j} \to H^k(L) \to E_\infty^{i,k-i} \to 0 \,.
$$ 
\end{lem}

\begin{proof}
Since, by hypothesis, we have $E_\infty^{r, k-r} = 0$ for all $r \in \{0, \dots, k\} \backslash \{ i, j \}$, it follows that in the filtration $(A^k_l)_{l = 0}^{k+1}$ of $H^k(L)$ we have $A^k_r = A^k_{r+1}$ for all $r \in \{0, \dots, k\} \backslash \{ i, j \}$ and, therefore,
$$
A^k_r = 
\begin{cases}
H^k(L), & \text{ if } r \in \{0, \dots, i \}, \\
E_\infty^{j, k-j}, & \text{ if } r \in \{i + 1, \dots, j \}, \\
0, & \text{ if } r \in \{j + 1, \dots, k \} \,.
\end{cases}
$$
Now, since $E_\infty^{i, k-i} \cong A^k_i/A^k_{i+1}$, we may conclude that there is a short exact sequence
$$
0 \to E_\infty^{j, k-j} \to H^k(L) \to E_\infty^{i,k-i} \to 0 \,,
$$ 
as desired.
\end{proof}

Since $M^n$ is a highly connected rational homology sphere, Lemma \ref{L:SES} can be used to extract a useful relationship between the cohomology groups of $L$, $F$ and $M^n$.

\begin{lem}
\label{L:edge_homom}
For every $j \in \{0, \dots, m-1\}$, there exists an isomorphism $H^j(L) \cong H^j(F)$, while $H^m(L) \cong E_\infty^{0,m} \cong \ker ( d_{m+1}^{0,m} ) \In H^m(F)$.  

Moreover, for every $k \geq 0$ there exists an injective homomorphism $ E_\infty^{k, 0} \to H^{k}(L) $.  In particular, when $k = m+1$ there exists a short exact sequence 
$$
0 \to E_\infty^{m+1, 0} \to H^{m+1}(L) \to E_\infty^{0,m+1} \to 0 \,,
$$ 
where the group $E_\infty^{m+1, 0}$ is finite and given by
$$
E_\infty^{m+1, 0} \cong H^{m+1}(M^n) / \im(d_{m+1}^{0,m}) \,.
$$
\end{lem}

\begin{proof}
Recall that the entries on the $E_2$-page of the spectral sequence are given by $E_2^{p,q} = H^p(M; H^q(F))$.  From the universal coefficient theorem, it thus follows that there is a short exact sequence
$$
0 \to \Ext(H_{p-1}(M), H^q(F)) \to E_2^{p,q} \to \Hom(H_p(M), H^q(F)) \to 0 \,.
$$
Now, since $M^n$ is an $(m-1)$-connected rational homology sphere, we observe that, for all $q \in \N$, we have $E_2^{p,q} = 0$, whenever $p \in \{1, \dots, m-1 \}$, and $E_2^{m,q} \cong \Hom(H_m(M), H^q(F))$.  Observe further that, since $H_m(M)$ is a finite group and $H^0(F)$ and $H^1(F)$ are always free abelian groups (see Table \ref{table:cohomF}), it must follow that $E_2^{m,q} = 0$ whenever $q \in \{0,1\}$.

Thus, for $j \in \{0, \dots, m\}$, the only possible non-trivial entry along each diagonal $\{ E_\infty^{p,q} \mid p + q = j \}$ is the term $E_\infty^{0,j}$, meaning that $H^j(L) \cong E_\infty^{0,j}$ for each $j \in \{0, \dots, m\}$.  However, when $j \in \{0, \dots, m-1\}$, all differentials originating from $E_*^{0,j}$ are necessarily trivial, which implies that $H^j(L) \cong E_\infty^{0,j} \cong  E_2^{0,j} = H^j(F)$ for all $j \in \{0, \dots, m-1\}$.  On the other hand,  the only possible non-trivial differential $d_r^{0, m} : E_r^{0,m} \to E_r^{r, m-r+1} $ involving a term of the form $E_*^{0,m}$ occurs when $r = m+1$, implying that 
$$
H^m(L) \cong E_\infty ^{0,m} \cong E_{m+2}^{0,m} = \ker( d_{m+1}^{0, m} ) \In E_{m+1}^{0,m} \cong E_2^{0,m} = H^m(F) \,.
$$

\begin{figure}
\centering
\scalebox{0.86}{
\begin{tikzpicture}
  \matrix (m) [ampersand replacement=\&, matrix of math nodes,
    nodes in empty cells,nodes={minimum width=5ex,
    minimum height=5ex,outer sep=-5pt},
    column sep=1ex,row sep=1ex]{
            F  \&    \&  *   \&   0  \& \cdots \& 0 \&  0 \&  \&  \&    \\
                \&   m   \&  H^m (F)  \&   0  \& \cdots  \& 0 \&  0 \&  \&  \&    \\  
                \&   \node[left] {m-1} ; \&  H^{m-1}(F)   \&   0  \& \cdots  \& 0 \&  0 \&   \&  \&    \\ 
                \&   \vdots   \&  \vdots   \&  \vdots   \& \vdots  \&  \vdots \& \vdots \&   \&   \&     \\
                \&   2   \&  H^2(F)   \&   0  \& \cdots  \& 0 \& 0 \& * \& * \&    \\   
                \&   1   \&  H^1(F)   \&   0  \& \cdots  \& 0 \& 0  \& 0  \& * \&      \\ 
                \&   0   \&  \Z   \&  0   \& \cdots  \& 0 \& 0 \& 0 \& H^{m+1}(M) \&      \\   
                \&  \strut  \quad  \&  0    \&  1   \&  \cdots  \&  m-2 \& m-1 \& m  \& m+1 \&  \strut  M^{2m+1} \\
                 };
\draw[thick, color=lightgray] (m-1-2.north east) -- (m-8-2.east) ; 
\draw[thick, color=lightgray] (m-8-2.north) -- (m-8-9.north east) ; 

\draw[thick, -stealth] (m-2-3.south east) -- (m-7-9.north west);  
\node[xshift=4mm, yshift=1mm] at (m-4-5.north east) {$d^{0,m}_{m+1}$}; 

\end{tikzpicture}
}
\caption{The only non-trivial differential involving $H^m(F)$ and $H^{m+1}(M)$ is $d^{0,m}_{m+1}$.}
\label{F:specseq}
\end{figure}

The fact that there is an injection $ E_\infty^{k, 0} \to H^{k}(L) $ for all $k \geq 0$ (called the edge homomorphism) is a trivial consequence of the definition of the filtration of $H^{k}(L)$, since $E_\infty^{k,0} \cong A^k_k \In H^k(L)$ for any $k \geq 0$.  

From the discussion above, the only possible non-trivial entries along the diagonal $\{ E_2^{p,q} \mid p + q = m+1 \}$ are $E_2^{0,m+1}$ and $E_2^{m+1, 0}$.  Therefore, the only possible non-trivial entries along the diagonal $\{ E_\infty^{p,q} \mid p + q = m+1 \}$ are $E_\infty^{0,m+1}$ and $E_\infty^{m+1, 0}$.  By Lemma \ref{L:SES}, this establishes the existence of the desired short exact sequence.

Finally, from the observations above, the only possible non-trivial differential $ d_{r}^{m+1-r, r-1} : E_{r}^{m+1-r, r-1} \to E_{r}^{m+1,0} $ involving a term of the form $E_*^{m+1,0}$ occurs when $r = m+1$.  Therefore, $E_{m+1}^{m+1,0} \cong E_{2}^{m+1,0} = H^{m+1}(M)$ and 
\begin{align*}
E_\infty^{m+1,0} \cong E_{m+2}^{m+1,0} &= E_{m+1}^{m+1,0} / \im(d_{m+1}^{0,m}) \\
&\cong H^{m+1}(M) / \im(d_{m+1}^{0,m}) \,.
\end{align*}
\end{proof}

It is now possible to make an initial observation of the strong topological restrictions on highly connected rational homology spheres which admit a double disk-bundle decomposition.

\begin{thm}
\label{T:max2}
Suppose that $M^n$ is an $(m-1)$-connected, rational homology sphere of dimension $n = 2m+1$,  $m \geq 2$.  If $M^n$ admits a double disk-bundle decomposition $DB_- \cup_L DB_+$ with $B_\pm$ connected, then a minimal generating set for the torsion group $H_{m}(M^n)$ consists of at most two elements.  In particular, if $H^{m}(L)$ has torsion, then $n \equiv 1$ mod $4$, $H^m(L) = \Z_2$ and $H_m(M^n) \in \{0, \Z_2 \}$.
\end{thm}

\begin{proof}
By Poincar\'e duality, it is sufficient to show that a minimal generating set for the group $H^{m+1}(M^n)$ consists of at most two elements.  To this end, recall, from Lemma \ref{L:edge_homom}, that there exists an injective homomorphism
\beq
\label{E:injection}
H^{m+1}(M) / \im(d_{m+1}^{0,m}) \cong E_\infty^{m+1, 0} \to H^{m+1}(L) \,.
\eeq

Suppose first that $H^{m+1}(L)$ is torsion free.  Then, since $H^{m+1}(M^n) \cong H_{m}(M^n)$ is a finite group, it follows that 
$$
H^{m+1}(M) / \im(d_{m+1}^{0,m}) = 0
$$ 
and, hence, that the differential
$$
d_{m+1}^{0, m} : E_{m+1}^{0, m} \to E_{m+1}^{m+1,0}  
$$
is surjective.  Now, while in the proof of Lemma \ref{L:edge_homom} it was established that $E_{m+1}^{0, m} = H^m(F)$ and $E_{m+1}^{m+1,0} = H^{m+1}(M)$.  Therefore, since it is known from Table \ref{table:cohomF} that a minimal generating set for each cohomology group of $F$ consists of at most two elements, the surjectivity of the differential implies that a minimal generating set for $H^{m+1}(M)$ must also consist of at most two elements.

Suppose, on the other hand, that $H^{m+1}(L)$ has torsion.  Then, since $\dim(L) = n-1 = 2m$, it follows from Poincar\'e duality and the universal coefficient theorem that the torsion subgroup of $H^m(L)$ is isomorphic to that of $H^{m+1}(L)$.  By Lemma \ref{L:edge_homom}, it follows that $H^m(F)$ has torsion.  Corollary \ref{C:torF} now implies that $n \equiv 1$ mod $4$ and that $H^m(F) = \Z_2$ (note, in particular, that the case $n=7$ cannot occur).

Therefore, since $H^m(L)$ has been assumed to be non-trivial, it may be concluded from Lemma \ref{L:edge_homom} that $H^m(L) \cong H^m(F) = \Z_2$ (since $H^m(L) \In H^m(F)$) and, hence, the torsion subgroup of $H^{m+1}(L)$ must also be isomorphic to $\Z_2$.  Moreover, the image of the differential $d_{m+1}^{0,m}$ is trivial, since $\im(d_{m+1}^{0,m}) \cong H^m(F)/H^m(L) = 0$.  Consequently, the injection \eqref{E:injection} yields an injection
$$
H^{m+1}(M) \to \Z_2 \In H^{m+1}(L) \,,
$$
as desired.
\end{proof}

By Theorem \ref{T:max2}, it is necessary to understand double disk-bundle decompositions where $H^m(L)$ is free abelian.

\begin{lem}
\label{L:Lfreeab}
If $H^m(L)$ is free abelian, then $H^{m+1}(M) \cong H^m(F)/ H^m(L)$ and either 
\begin{enumerate}
\item $H^m(L) = H^m(B_\pm)= 0$ and $H^m(F) \in \{0, \Z_2\}$, or

\vspace{2mm}
\item $H^m(L) \cong \Z^2$, $H^m(B_\pm) \cong \Z$ and $H^m(F) \cong \Z^2$.
\end{enumerate}
In particular, if $H^m(L) = 0$ then $H^{m+1}(M) \in \{0, \Z_2\}$.
\end{lem}

\begin{proof}
By Poincar\'e duality and the universal coefficient theorem, $H^m(L)$ being free abelian means that the same must be true for $H^{m+1}(L)$.  Thus, as discussed in the proof of Theorem \ref{T:max2}, there is an isomorphism
\beq
\label{E:iso}
H^{m+1}(M) = \im(d_{m+1}^{0,m}) \cong H^m(F)/ H^m(L) \,.
\eeq
Since $H^{m+1}(M)$ is finite, this implies that $b_m(L) = b_m(F)$, i.e.~that the $m^\text{th}$ Betti numbers of $F$ and $L$ must be equal.  Moreover, Table \ref{table:cohomF} then yields the inequality $b_m(L) = b_m(F) \leq 2$, with equality if and only if $H^m(F) \cong \Z^2$.

Since $M^n$ is a rational homology sphere, the commutativity of the braid diagram \eqref{E:braid} yields rational short exact sequences
$$
0 \to H^j(B_\pm; \Q) \to H^j(L; \Q) \to H^{j+1}(M, B_\mp; \Q) \to 0
$$
for all $j \in \{1, \dots, 2m-1\}$, from which it follows via the isomorphisms \eqref{E:isoms} that the Betti numbers of $L$ and $B_\pm$ satisfy $b_j(L) = b_j(B_\pm) + b_{2m-j}(B_\pm)$ for all $j \in \{1, \dots, 2m-1\}$.  In particular, setting $j = m$ implies that $b_m(L) = 2 \, b_m(B_\pm)$, which, in combination with the inequality $b_m(L) = b_m(F) \leq 2$ above, yields that the only possibilities for the (free-abelian, integral) group $H^m(L)$ are $H^m(L) = 0$ and $H^m(L) \cong \Z^2$. 

Suppose first that $H^m(L) = 0$.  Then $H^m(F) \cong H^{m+1}(M)$ by \eqref{E:iso}, meaning that $H^m(F)$ is finite, while $H^m(B_\pm)$ must also be finite, since $0 = b_m(L) =2 \, b_m(B_\pm)$.  Hence, as in the proof of Theorem \ref{T:max2}, the only possibilities for $H^m(F)$ are $H^m(F) = 0$ or $H^m(F) = \Z_2$ (in which case $n = 2m+1 \equiv 1$ mod $4$). 

Suppose, on the other hand, that $H^m(L) \cong \Z^2$.  By the above discussion, this implies that $H^m(F) \cong \Z^2$ and that $b_m(B_\pm) = 1$.

Therefore, it remains only to show that $H^m(B_\pm)$ is free abelian in all cases.  However, this follows easily from the commutativity of the braid diagram \eqref{E:braid}, since $M^n$ being $(m-1)$-connected implies that $H^{m}(M) = 0$, which in turn leads to portions
$$
0 \to H^m(B_\pm) \to H^m(L) 
\to \cdots
$$
in the exact sequences for the pairs $(DB_\pm, L)$; that is, the groups $H^m(B_\pm)$ inject into the free abelian group $H^m(L)$.
\end{proof}

All the elements are now in place to prove Theorem \ref{thm:main2}.

\begin{thm}
\label{T:cyclic}
Suppose that $M^n$ is an $(m-1)$-connected, rational homology sphere of dimension $n = 2m+1$,  $m \geq 2$.  If $M^n$ admits a double disk-bundle decomposition $DB_- \cup_L DB_+$ with $B_\pm$ connected, then $H_{m}(M^n)$ is cyclic.
\end{thm}

\begin{proof}
By Theorem \ref{T:max2} and Lemma \ref{L:Lfreeab}, the only situation in which the finite group $H_m(M^n) \cong H^{m+1}(M^n)$ could possibly fail to be cyclic is when $H^m(L) \cong \Z^2$.  Assume, therefore, that $H^m(L) \cong \Z^2$ and consider the braid diagram \eqref{E:braid} when $j = m$.  Since $M^n$ is $(m-1)$-connected, it is clear that $H^m(M) = 0$ and, by Lemma \ref{L:Lfreeab}, the assumption that $H^m(L) \cong \Z^2$ yields that $H^m(B_\pm) = \Z$.  The universal coefficient theorem and the isomorphisms \eqref{E:isoms} then ensure that $H^{m+1}(M, B_\pm) \cong H_m(B_\mp)$ is of rank one, with torsion isomorphic to that of $H^{m+1}(B_\mp)$.  To prove the theorem, it now suffices to show that the groups $H^{m+1}(B_\pm)$ are torsion free, since this would yield short exact sequences 
$$
0 \to H^m(B_\pm) \cong \Z \to H^{m+1}(M, B_\pm) \cong \Z \to H^{m+1}(M) \to 0
$$
in the braid diagram \eqref{E:braid}, as desired.

Thus, it remains to show that the groups $H^{m+1}(B_\pm)$ are torsion free whenever $H^m(L) \cong \Z^2$.  To that end, recall from Lemma \ref{L:Lfreeab} that $H^m(F) \cong \Z^2$.  Therefore, by Corollary \ref{C:torF} and Table \ref{table:cohomF}, there are three cases to consider:
\begin{enumerate}
\item 
\label{i:same+or}
$\ell_+ = \ell_-$, with $\sph^{\ell_\pm} \to L \to B_\pm$ both orientable and 
$$
H^j(F) \cong \begin{cases}
\Z, & j = 0 \,, \\
\Z^2, & j > 0 \ \text{ and } \ j \equiv 0 \!\! \mod \ell_+ \,, \\
0, & \text{otherwise};
\end{cases}
$$

\item 
\label{i:diff+or}
$\ell_+ \neq \ell_-$, with $\sph^{\ell_\pm} \to L \to B_\pm$ both orientable and 
$$
H^j(F) \cong \begin{cases}
\Z, & j = 0 \ \text{ or } \ j \in \{\ell_\pm\} \!\! \mod \ell_+ + \ell_- \,, \\
\Z^2, & j > 0 \ \text{ and } \ j \equiv 0 \!\! \mod \ell_+ + \ell_- \,, \\
0, & \text{otherwise};
\end{cases}
$$

\item 
\label{i:non-or}
$\ell_+ = \ell_- = 1$ and $n = 7$, with $\sph^{\ell_\pm} \to L \to B_\pm$ both non-orientable and 
$$
H^j(F) \cong \begin{cases}
\Z, & j = 0\,, \\
\Z^2, & j > 0 \ \text{ and } \ j \equiv 0 \!\! \mod 3\,, \\
0, & j \equiv 1 \!\! \mod 3\,, \\
\Z_2^2, & j \equiv 2 \!\! \mod 3\,.
\end{cases}
$$
\end{enumerate}

Furthermore, recall from Lemma \ref{L:edge_homom} that $H^j(L) \cong H^j(F)$ for all $j \in \{0, \dots, m-1\}$.  In particular, in cases \eqref{i:same+or} and \eqref{i:diff+or} one easily deduces from Poincar\'e duality and the universal coefficient theorem that the cohomology groups of $L$ are free abelian and given by
\begin{equation}
\label{E:Lgps}
H^j(L) \cong \begin{cases}
H^j(F), & j \in \{0, \dots, m-1\}\,, \\
\Z^2, & j =m \,, \\
H^{2m-j}(F), & j \in \{m+1, \dots, 2m\} \,.
\end{cases}
\end{equation} 

\noindent
\underline{Case \eqref{i:same+or}:} 
By Corollary \ref{C:dims} and since $n$ is odd, the homotopy fiber $F$ has a loop-space factor $\Omega \sph^k$, with $n = k$, if $k$ is odd, and $n = 2k-1$, if $k$ is even.  Moreover, since the bundles $\sph^{\ell_\pm} \to L \to B_\pm$ are both orientable and since $n = 2m+1 \geq 5$ by hypothesis, it follows from Table \ref{table:Qtype} that $\ell_\pm \geq 2$ and that either $k = \ell_\pm + 1$ or $k = 2 \ell_\pm + 1$ or
$$
(\ell_\pm, k) \in \{(2,7), (2,9), (2,13), (4,13), (4,17), (4,25), (8, 25)\}.
$$
Now, since $H^m(F) \cong \Z^2$, there exists some $r \in \N$ such that $m = r \ell_\pm$.  However, since $n = k$, if $k$ is odd, and $n = 2k-1$, if $k$ is even, the above-listed possibilities for $k$ yield that $(m,k) = (\ell_\pm, \ell_\pm + 1)$ (with $m = \ell_\pm$ odd) or $(m,k) = (\ell_\pm, 2\ell_\pm + 1)$ or $m = r \ell_\pm$ with $r \in \{2,3\}$ and
$$
(\ell_\pm, k) \in \{(2,9), (2,13), (4,17), (4,25)\}.
$$
Observe, furthermore, that 
$$
\dim(B_\pm) = 2m - \ell_\pm = m + (r-1)\ell_\pm \,.
$$
Thus, it may be assumed that $m = r \ell_\pm$ with $r \in \{2,3\}$, since $r=1$ gives $\dim(B_\pm) = m$ and, hence, that $H^{m+1}(B_\pm) = 0$.

Consider now the Gysin sequences 
$$
\cdots \to H^{j+\ell_\pm}(L) \to H^j(B_\pm) \to H^{j+\ell_\pm + 1}(B_\pm) \to H^{j+\ell_\pm + 1}(L) \to \cdots
$$
associated to the sphere bundles $\sph^{\ell_\pm} \to L \to B_\pm$.  If $j=m+1$ then, by \eqref{E:Lgps} and since $\ell_\pm \geq 2$, 
$$
H^{m+\ell_\pm +1}(L) \cong H^{m-\ell_\pm -1}(F) = H^{(r-1)\ell_\pm - 1}(F) = 0 
$$
because $(r-1)\ell_\pm - 1 \not\equiv 0$ mod $\ell_+$.  Thus, there is an injection $H^{m+1}(B_\pm) \to H^{m+\ell_\pm+2}(B_\pm)$.  However, since $m+\ell_\pm+2 = \dim(B_\pm) + (2-r)\ell_\pm+2$, it follows that $H^{m+\ell_\pm+2}(B_\pm) = 0$, if $r = 2$, and $H^{m+\ell_\pm+2}(B_\pm) = \Z$, if $r=3$ and $\ell_\pm = 2$.  In either case, it follows that $H^{m+1}(B_\pm)$ must be free abelian.  

Finally, let $r = 3$ and $\ell_\pm > 2$ (i.e.~$\ell_\pm  = 4$).  If $j = m + \ell_\pm + 2$ in the Gysin sequence above, then 
$$
H^{m+2\ell_\pm +2}(L) \cong H^{m-2\ell_\pm -2}(F) = H^{\ell_\pm - 2}(F) = 0
$$
by \eqref{E:Lgps}, while $H^{m + 2\ell_\pm + 3}(B_\pm) = 0$ because $m + 2\ell_\pm + 3 = \dim(B_\pm) + 3$.  Thus, it must again follow that $H^{m+\ell_\pm+2}(B_\pm) = 0$ and, via the above-mentioned injection, that $H^{m+1}(B_\pm) = 0$.

\vspace{2mm}
\noindent
\underline{Case \eqref{i:diff+or}:} 
As in case \eqref{i:same+or}, Corollary \ref{C:dims} yields that the homotopy fiber $F$ has a loop-space factor $\Omega \sph^k$, with $n = k$, if $k$ is odd, and $n = 2k-1$, if $k$ is even.  Assume now, without loss of generality, that $\ell_- < \ell_+$.  Then, since the bundles $\sph^{\ell_\pm} \to L \to B_\pm$ are both orientable, it follows from Table \ref{table:Qtype} that $1 \leq \ell_- < \ell_+ $ and that $k = \ell_- + \ell_+ + 1$.

Observe that $k$ must be even since, otherwise, we have $2 \leq m < n-1 = k-1 = \ell_- + \ell_+$, meaning that $H^m(F) \not\cong \Z^2$, in contradiction to Lemma \ref{L:Lfreeab}.  Therefore, we may assume that $k$ is even and that $n = 2k-1$, yielding in turn that $m = k-1 = \ell_- + \ell_+$ is odd.  In particular, it follows that $\ell_- + 1 < \ell_+ + 1 < m+1$.  Therefore, choosing $j = \ell_\pm + 1$, respectively, in the braid diagram \eqref{E:braid} yields injections $H^{\ell_\pm + 1}(B_\pm) \to H^{\ell_\pm + 1}(L)$ induced by the inclusion maps $L \to DB_\pm$.  Since these inclusion maps are homotopic to the bundle projections $L \to B_\pm$, it follows that the corresponding induced homomorphisms $H^{\ell_\pm + 1}(B_\pm) \to H^{\ell_\pm + 1}(L)$ in the Gysin sequences for the bundles $\sph^{\ell_\pm} \to L \to B_\pm$ are also injective.  Hence, the homomorphisms $\Z \cong H^0(B_\pm) \to H^{\ell_\pm + 1}(B_\pm)$ given by cupping with the respective Euler classes must be trivial.  In other words, the Euler classes themselves are trivial and the Gysin sequences split into short exact sequences.  In particular, we obtain short exact sequences
$$
0 \to H^{m + 1}(B_\pm) \to H^{m + 1}(L) \to H^{\ell_\mp + 1}(B_\pm) \to 0
$$
and, since $H^{m+1}(L)$ is free abelian, it follows that $H^{m+1}(B_\pm)$ is also free abelian, as desired.

\vspace{2mm}
\noindent
\underline{Case \eqref{i:non-or}:} 
In this case, we have $m = 3$ and that both $B_\pm$ are non-orientable.  It will be shown that $H^4(B_\pm) = 0$.  

Given that $H^j(L) \cong H^j(F)$ for all $j \in \{0, \dots, 2\}$, by Lemma \ref{L:edge_homom}, it follows from Poincar\'e duality and the universal coefficient theorem that $H^4(L) =0$ and $H^5(L) \cong \Z_2^2$.  
Now, since $M$ is a $2$-connected rational homology sphere, it is clear from considering the braid diagram \eqref{E:braid} with rational coefficients that $H^2(B_\pm)$ and $H^4(B_\pm)$ must both be finite and, hence, isomorphic to the torsion subgroups of $H_1(B_\pm)$ and $H_3(B_\pm)$, respectively, by the universal coefficient theorem.  Moreover, since $H^3(B_\pm) \cong \Z$ (by Lemma \ref{L:Lfreeab}), it follows that $H_2(B_\pm) = 0$.

Since $M$ is $2$-connected, the long exact homotopy sequence for $F \to L \to M$ and Table \ref{table:Qtype} together yield $\pi_1(L) \cong Q_8$.  On the other hand, if $\iota_\pm : \sph^1 \to L$ denote the respective inclusions of a fiber in the bundles $\sph^1 \to L \to B_\pm$, then Proposition \ref{prop:pi1} ensures that $\pi_1(L)$ is generated by the non-trivial cyclic subgroups $(\iota_\pm)_*(\pi_1(\sph^1))$.  Thus,  the long exact homotopy sequences for the bundles $\sph^1 \to L \to B_\pm$ yield that $0 \neq \pi_1(B_\pm) \cong \pi_1(L) / (\iota_\pm)_*(\pi_1(\sph^1))$ and, hence, $\pi_1(B_\pm)$ must be isomorphic to the images of the cyclic subgroups $(\iota_\mp)_*(\pi_1(\sph^1))$ under the respective quotient homomorphisms.  This implies that $\pi_1(B_\pm) \cong \Z_2$, the only possible non-trivial, cyclic quotient of $Q_8$.  In particular, it follows that $H^2(B_\pm) \cong H_1(B_\pm) \cong \pi_1(B_\pm) \cong\Z_2$.

Let $\wt B_\pm$ be the (five-dimensional) universal covers of $B_\pm$, respectively.  Poincar\'e duality implies that $H^4(\wt B_\pm) \cong H_1(\wt B_\pm) = 0$ and, hence, the homomorphisms $H^4(B_\pm) \to H^4(\wt B_\pm)$ induced by the respective two-fold covering maps are trivial.  By considering the transfer homomorphism associated to the coverings, it now follows that the groups $H^4(B_\pm)$ consist entirely of $2$-torsion; that is, $H^4(B_\pm) \cong \Z_2^{r_\pm}$, for some $r_\pm \in \N \cup \{0\}$.  Furthermore, applying the universal coefficient theorem, observe that $H_3(B_\pm) \cong \Z \oplus \Z_2^{r_\pm}$.

Consider now cohomology with coefficients in $\Z_2$.  By applying the Universal Coefficent Theorem, we find that $H^j(L; \Z_2) \cong \Z_2^2$ for all $j \in \{1, \dots, 5\}$ and that
$$
H^j(B_\pm; \Z_2) \cong 
\begin{cases} 
\Z_2\,, & j = 0, 1, 2, \\
\Z_2^{1 + r_\pm}\,, & j = 3\,.
\end{cases}
$$
Therefore, the portion 
$$
0 \to H^1(B_\pm; \Z_2) \to H^1(L; \Z_2) \to H^0(B_\pm; \Z_2) \xrightarrow{\smile e_\pm} H^2(B_\pm; \Z_2) \to \cdots 
$$
of the Gysin sequence for $\sph^1 \to L \to B_\pm$ yields that the respective $\Z_2$-Euler classes $e_\pm \in H^2(B_\pm; \Z_2)$ must be trivial and, hence, that the Gysin sequence splits into a family of short exact sequences
$$
0 \to H^j(B_\pm; \Z_2) \to H^j(L; \Z_2) \to H^{j-1}(B_\pm; \Z_2) \to 0\,.
$$
By setting $j = 3$, we obtain $H^3(B_\pm; \Z_2) \cong \Z_2$, which forces $r_\pm = 0$ and, consequently, we arrive at $H^4(B_\pm) = 0$, as desired.
\end{proof}


\end{document}